\documentclass{amsart}
\usepackage{graphicx}
\usepackage{BUSSPROOFS}
\usepackage{amssymb}
\usepackage{amsfonts}
\usepackage{latexsym}
\usepackage{amstext}
\usepackage{amsmath}
\usepackage{amsthm}
\usepackage{xcolor}
\usepackage{newlfont}
\usepackage{enumitem}
\newtheorem{thm}{Theorem}[section]
\newtheorem{cor}[thm]{Corollary}
\newtheorem{lem}[thm]{Lemma}

\theoremstyle{definition}
\newtheorem{defn}[thm]{Definition}
\theoremstyle{remark}
\newtheorem{rem}[thm]{Remark}

\numberwithin{equation}{section}
\newenvironment{proofof}[1]{\begin{trivlist}\item[\hskip\labelsep{\sc
Proof~of~{#1}.\ }]}{\hspace*{\fill} {\sc qed}\end{trivlist}}

\renewcommand{\a}{\alpha}
\newcommand{\be}{\beta}
\newcommand{\g}{\gamma}
\newcommand{\de}{\delta}

\begin{document}

\title{Uniform Interpolation in provability logics}%
\author{Marta B\'{i}lkov\'{a}}%
\address{Institute of Computer Science, CAS in Prague}%
\email{bilkova@cs.cas.cz}%

\thanks{The work of this paper has been supported by the project No. P202/11/1632 of the Czech Science Foundation. I would like to thank Nick Bezhanishvilli for pointing up the topic years ago, Albert Visser for inspiration and discussions on the topic, and Rosalie Iemhoff and Tadeusz Litak for encouraging me to write this up once again.}

\dedicatory{To Albert Visser}
\begin{abstract}
We prove the uniform interpolation theorem in modal provability logics \textbf{GL} and \textbf{Grz} by a proof-theoretical method, using analytical and terminating sequent calculi for the logics. The calculus for G\"{o}del-L\"{o}b's logic \textbf{GL} is a variant of the standard sequent calculus of \cite{SV82}, in the case of Grzegorczyk's logic \textbf{Grz}, the calculus implements an explicit loop-preventing mechanism inspired by work of Heuerding \cite{Heu96,Heu98}.  
\end{abstract}
\maketitle
\EnableBpAbbreviations
\section{Introduction}

\subsubsection{Uniform interpolation}

Uniform Interpolation Property for a logic $L$ is a strong interpolation property, stating that, for any formula $\alpha$ and any propositional variable $p$, there is a post-interpolant $\exists p(\alpha)$ not containing $p$, such that ${\vdash_L\alpha\to\exists p(\alpha)}$, and 
$\vdash_L\alpha\to\beta$ implies $\vdash_L\exists p(\alpha)\to\beta$ for each $\beta$ not containing $p$.
Similarly, for any $\beta$ and $p$ there is a pre-interpolant $\forall p(\beta)$ not containing $p$, such that $\vdash_L \forall p(\beta)\to\beta$ and 
$\vdash_L\alpha\to\beta$ implies $\vdash_L \alpha\to \forall p(\beta)$ for each $\alpha$ not containing $p$.
Uniform interpolation property entails Craig interpolation property, and uniform interpolants are unique up to the~provable equivalence, they are the the~minimal and the~maximal interpolants of a given implication w.r.t. the~provability ordering.

While for classical propositional logic, and also for other locally tabular logics like modal logic \textbf{S5}, uniform interpolation property is easily obtained, in other logics it is not the case. Interest in the topic arose with a seminal work by Pitts \cite{Pit92}, who proved uniform interpolation for intuitionistic propositional logic using a terminating sequent calculus.
For modal logic \textbf{K} uniform interpolation was first proved by Visser \cite{Viss96} and Ghilardi \cite{Ghil95}, for provability logic \textbf{GL} by Shavrukov \cite{Shav93}. The~failure of uniform interpolation in modal logic
\textbf{S4}, which applies to \textbf{K4} as well, was proved by Ghilardi and Zawadowski \cite{GZ95b}. 
More recent is a proof for monotone modal logic by Venema and Santocanale \cite{VenSan10}, using a coalgebraic perspective: uniform interpolants are constructed via erasing variable in a disjunctive normal form. 
This relates to the way the problem of computing uniform interpolants is understood in Artificial intelligence, which is variable forgetting. Similar motivation, but different approach based on resolution calculi and conjunctive normal forms, is applied to modal logic \textbf{K }by Herzig et al. \cite{Herz08}.

As the notation suggests, uniform interpolants relate to a certain type of propositional quantifiers:  if propositional quantifiers satisfying at least the usual quantifier axioms and rules are expressible in the language, they are the uniform interpolants. On the other hand, if we construct uniform interpolants so that the construction commutes with substitutions, we can use them to interpret the propositional quantifiers, precisely as was done by Pitts' in \cite{Pit92}.
Visser \cite{Vis96} proved uniform interpolation for various modal logics, provability logics \textbf{GL} and \textbf{Grz} among them, via a semantical argument which yields a~semantic
characterization of the~resulting quantifiers as \emph{bisimulation}
quantifiers: from the~semantic point of view, quantifying over $p$, they quantify over possible valuations of $p$ in models bisimilar to the current one up to $p$. A~complexity bound of uniform interpolants in terms of $\Box$-depth is obtained in the~proof, however, the~proof does not provide us with a~direct construction of the~interpolants. 

Using a method similar to Pitts' and using sequent calculi for modal logics, the author proved effective uniform interpolation for modal logics \textbf{K }and \textbf{T} in \cite{Bil06,Bilt06}. The thesis \cite{Bilt06} also contains proofs of uniform interpolation for provability logics \textbf{GL} and \textbf{Grz} which are reconsidered in this paper. The reason we came back to the topic is a
recent interest in uniform interpolation in modal and modal intuitionistic logic by Iemhoff \cite{Iemh15}.  

\subsubsection{Provability modal logics}

In this paper, we concentrate solely on the G\"odel-L\"ob's provability logic \textbf{GL}, and Grzegorczyk!s logic \textbf{Grz}, also known as \textbf{S4Grz}. The~main reference for provability modal logics and their properties is Boolos' book \cite{Boo93}, for a~history of provability logic see also \cite{SB91}.

The logic \textbf{GL} is a normal modal logic, extending the basic modal logic \textbf{K }with the L\"ob's axiom 
$$
L:\ \Box(\Box p\to p)\to \Box p.
$$
It is known to be complete with respect to transitive and conversely well-founded Kripke frames.
The logic \textbf{Grz} is a normal modal logic, extending the basic modal logic \textbf{K }with the axiom $T:\ \Box p\to p$, and the Grzegorczyk's axiom 
$$
Grz:\ \Box(\Box(p\to\Box p)\to p)\to \Box p.
$$
It is known to be complete with respect to transitive, reflexive and conversely well-founded Kripke frames. Both logics have the finite model property as well and are therefore decidable, as was shown e.g. in \cite{Avr84}.

The~$\Box$ modality of G\"{o}del-L\"{o}b's logic \textbf{GL} can be interpreted as
formalized provability in an~arithmetical recursively axiomatizable theory $T$: assume an axiomatization of $T$ is expressed by a~sentence $\tau$ and consider a~standard proof predicate $Pr_{\tau}(\bar{\varphi})$ for $T$. An~arithmetical interpretation of modal formulas is a~function from propositional variables to arithmetical sentences such that it
commutes with logical connectives, and $e(\bot)=(0=S(0))$, and
$e(\Box\a)=Pr_{\tau}(\overline{e(\a)})$.
Arithmetical completeness is established as the following statement:
$$
\vdash_{H_{GL}}\a\quad\mbox{iff}\quad \forall e (T\vdash e(\a)). 
$$
G\"{o}del-L\"{o}b's logic \textbf{GL} was proved to be arithmetically complete for
Peano arithmetic by Solovay \cite{Sol76}. Later it was shown that it is the~logic of provability of a~large family of reasonable formal theories.

Using the above interpretation of \textbf{GL}, we obtain the~following arithmetical
interpretation of Grzegorczyk's logic: an~arithmetical
interpretation of modal formulas is as before, only now $e(\Box\a)=Pr_{\tau}(\overline{e(\a)})\wedge\overline{e(\a)}.$

\section{Calculi}

To prove the uniform interpolation theorem we use sequent calculi with good structural properties.
The~particular form of sequent calculi has been chosen for proof-search related manipulations. In particular, we use finite multisets of formulas to formulate a~sequent, a notation which does not hide contractions (contraction rules are not part of the definition and are to be proved admissible rules), we use a~definition without the~cut rule (which is to be proved admissible), and structural rules of contraction and weakening are built in logical rules and axioms.
Since the~proof of the uniform interpolation theorem contained in the~next
section is closely related to termination of a~proof-search in the~calculi, we will devote some space in this section to explain proof-search in provability logics and its termination. Namely, we employ simple implicit loop-preventing mechanisms provided naturally by diagonal formulas, and in the case of Grzegorczyk's logics also an explicit syntactic loop-preventing mechanism to avoid reflexive loops due to the presence of the $T$ axiom.

We assume the reader is familiar with basics on sequent calculi as contained e.g. in the Schwichtenberg's and Troelstra's book
\cite{ST96}. For sequent calculi of modal logics having arithmetical interpretation we refer to Sambin and Valentini's paper \cite{SV82}, or Avron's paper \cite{Avr84}.

\subsection{Preliminaries}

Formulas are given by the following grammar of the basic modal language, where atoms $p$ are taken from a fixed countable set of propositional variables:
$$
\alpha:= \bot\ |\ p\ |\ \neg\alpha\ |\ \alpha\wedge\alpha\ |\ \alpha\vee\alpha\ |\ \Box\alpha,
$$
the notion of subformulas is standard, and with the term atomic formula we refer to atoms as well as the constant $\bot$. We moreover define $\top=\neg\bot$ and $\alpha\to\beta = \neg\alpha\vee\beta$, $\diamondsuit\a=\neg\Box\neg\a$, and we use $\bigvee\emptyset= \bot$ and $\bigwedge\emptyset= \top$.  The \textit{weight} $w(\alpha)$ of a formula $\alpha$ is the number of symbols it contains, and the \textit{box-depth} $d(\alpha)$ of a formula $\alpha$ is defined as the maximum number of boxes along a branch in the corresponding formula tree.
By capital Greeks we denote finite multisets of formulas. Formally, $\Gamma$ is a function from the set of formulas to natural numbers with finite support (finitely many non-zero values), but we mostly use a relaxed notation and treat multisets as sets with multiple occurrences, in particular by $\alpha\in\Gamma$ we mean that $\Gamma(\alpha)> 0$. For a multiset $\Gamma$, we denote the underlying set by $\Gamma^{\circ}$. $\Box\Gamma$ denotes the multiset resulting from prefixing elements of $\Gamma$ with box while keeping the multiplicities intact. A \textit{sequent} is a syntactic object of the form $\Gamma\Rightarrow\Delta$, or, in the case of Grzegorczyk's logic, of the form $\Box\Sigma|\Gamma\Rightarrow\Delta$. The weight $w(\Gamma)$ of a multiset $\Gamma$ is the sum of the weights of elements in $\Gamma$, the weight $w(\Gamma,\Delta)$ of a sequent $\Gamma\Rightarrow\Delta$ is the sum $w(\Gamma)+w(\Delta)$.

A rule consists of a finite set of sequents called premises and a single sequent called the conclusion, rules with zero premises are called axioms. A calculus is given by a set of rule-schemes, a proof in the calculus is then a finite rooted tree labeled with sequents in such a way that leaves are labeled with axioms and labels of parent-children nodes respect correct instances of the rules of the calculus. The height of a proof is the height of the tree. A sequent is provable if there is a proof whose root is labeled with the sequent. 

We call a rule \textit{invertible} if whenever the conclusion is provable, then all its premises are provable as well, we call a rule height-preserving invertible if moreover the premises have proofs of at most the height of the proof of the conclusion. We call a rule \textit{admissible} if whenever its premises are provable, so is the conclusion, and height-preserving admissible if moreover the conclusion has a proof of at most the maximum of the heights of the proofs of the premises.

By a \textit{proof-search} in a calculus we mean a procedure based on applying rules of the calculus backwards to a sequent in such a way that, for a provable sequent, the resulting tree contains a proof of the sequent. We call a proof-search \textit{terminating} if it results in a finite tree. 
Particular instances of proof search will be defined later.

\subsection{Sequent calculus for \textbf{GL}}

The following calculus is a variant of the sequent calculus introduced in \cite{SV82}, and reconsidered in \cite{GR08,GR12} using multisets in place of sets.

\begin{defn}
Sequent calculus $G_{GL}$:

\begin{prooftree}
\AxiomC{$\Gamma,p\Rightarrow p,\Delta$}
\AxiomC{$\Gamma,\bot\Rightarrow\Delta$}
\noLine\BIC{ }
\end{prooftree}

\begin{prooftree}
\AxiomC{$\Gamma,\a,\be\Rightarrow\Delta$}
\RightLabel{$\wedge$-l}
\UnaryInfC{$\Gamma,\a\wedge\be\Rightarrow\Delta$}
\AxiomC{$\Gamma\Rightarrow\a,\be,\Delta$}
\RightLabel{$\vee $-r}
\UnaryInfC{$\Gamma\Rightarrow\a\vee\be,\Delta$}
\noLine\BinaryInfC{  }
\end{prooftree}

\begin{prooftree}
\AxiomC{$\Gamma,\a\Rightarrow\Delta$}
\RightLabel{$\neg$-r}
\UnaryInfC{$\Gamma\Rightarrow\neg\a,\Delta$}
\AxiomC{$\Gamma\Rightarrow\a,\Delta$}
\RightLabel{$\neg$-l}
\UnaryInfC{$ \Gamma,\neg\a\Rightarrow\Delta$}
\noLine\BinaryInfC{  }
\end{prooftree}

\begin{prooftree}
\AxiomC{$\Gamma\Rightarrow \a,\Delta$}\AxiomC{$\Gamma\Rightarrow\be,\Delta$}
\RightLabel{$\wedge$-r}
\BinaryInfC{$\Gamma\Rightarrow\a\wedge\be,\Delta$}
\AxiomC{$\Gamma,\a\Rightarrow\Delta$}
\AxiomC{$\Gamma,\be\Rightarrow \Delta$}
\RightLabel{$\vee$-l}
\BinaryInfC{$\Gamma,\a\vee \be\Rightarrow\Delta$}
\noLine\BinaryInfC{  }
\end{prooftree}

\begin{prooftree}
\AxiomC{$\Box\Gamma,\Gamma,\Box\a\Rightarrow\a$}
\RightLabel{$\Box_{GL}$}
\UnaryInfC{$\Box\Gamma,\Pi\Rightarrow\Box\a,\Delta$}
\end{prooftree}
In the~$\Box_{GL}$ rule, $\Pi$ contains only propositional
variables and $\Delta$ contains only propositional variables and
boxed formulas, and we call formulas $\Box\a$ as well as $\Box\Gamma$ principal formulas (formula occurrences). In the case of axioms, $p$ (resp. $\bot$) are principal formulas, and in the remaining rules, the principal formula is the one to which a connective is introduced.
\end{defn}

The propositional (non-modal) part of the calculus is a slight variant of the propositional part of the calculus \textbf{G3c} from \cite{ST96}. The propositional rules of the calculus are height-preserving invertible, for a proof of this fact we refer to \cite{ST96}. It is also not hard to prove that sequents of the form $\Gamma,\a\Rightarrow \a,\Delta$ are provable  for arbitrary $\a$.
\begin{lem}\label{lem:strGL} Weakening and contraction rules are
height-preserving admissible in $G_{GL}$.
\end{lem}
\begin{proofof}{Lemma \ref{lem:strGL}}

Weakening and contraction rules are:
\begin{prooftree}
\AXC{$\Gamma\Rightarrow\Delta$}
\RL{w-l}
\UIC{$\a,\Gamma\Rightarrow\Delta$}
\AXC{$\Gamma\Rightarrow\Delta$}
\RL{w-r}
\UIC{$\Gamma\Rightarrow\Delta,\a$}
\noLine\BIC{ }
\AXC{$\a,\a,\Gamma\Rightarrow\Delta$}
\RL{c-l}
\UIC{$\a,\Gamma\Rightarrow\Delta$}
\AXC{$\Gamma\Rightarrow\Delta,\a,\a$}
\RL{c-r}
\UIC{$\Gamma\Rightarrow\Delta,\a$}
\noLine\BIC{ }
\noLine\BIC{ }
\end{prooftree}
Proof is by induction on the weight of the principal formula $\a$ of the weakening (resp. contraction) inference, and for each weight on the height of the proof of the premise of the rule. We prove the admissibility of both the left and right weakening rules simultaneously, and the same applies to the two contraction rules.

\medskip\par\noindent
\textbf{Weakening:} For an~atomic weakening formula the~proof is obvious - note that
atomic weakening is built in axioms as well as in the~$\Box_{GL}$-rule. (The case of weakening-r by $\Box\a$ when the last inference is a~$\Box_{GL}$-inference is then also obvious since it is built-in the~rule as well.)
For a non atomic and not boxed formula we consider its main connective and use height-preserving
invertibility of the~corresponding propositional rule, then weakening by subformula(s) of
lower weight admissible by the induction hypothesis, and finally apply the~propositional rule. 

Let us therefore only spell out the step for a formula of the form $\a=\Box\be$. If it is not principal in the last step of the proof of the premise of the weakening, we simply permute the weakening upwards and use the induction hypothesis. So consider $\Box\be$ being the principle formula of a $\Box_{GL}$ inference. Notice that weakening-right by boxed formulas is built in the~$\Box_{GL}$ rule, therefore it is enough to consider weakening-left rule.

\medskip
\begin{itemize}
\item[-] Weakening-l by $\Box\be$, the~last step is a~$\Box_{GL}$ inference --- we permute the proof as follows:

\begin{prooftree}
\AXC{$\Box\Gamma,\Gamma,\Box\g\Rightarrow\g$}
\RightLabel{$\Box_{GL}\ \Longrightarrow$}
\UnaryInfC{$\Box\Gamma,\Pi\Rightarrow\Box\g,\Lambda$}
\RightLabel{w-l}
\UIC{$\Box\be,\Box\Gamma,\Pi\Rightarrow\Box\g,\Lambda$}
\AxiomC{$\Box\Gamma,\Gamma,\Box\g\Rightarrow\g$}
\RightLabel{w-l, i.h.}
\UnaryInfC{$\Box\be,\be,\Box\Gamma,\Gamma,\Box\g\Rightarrow \g$}
\RightLabel{$\Box_{GL}$}
\UnaryInfC{$\Box\be,\Box\Gamma,\Pi\Rightarrow\Box\g,\Lambda$}
\noLine\BinaryInfC{ }
\end{prooftree}
\end{itemize}

\par\noindent
\textbf{Contraction:} For $\a$ atomic, if the~premise is an~axiom,
the~conclusion is an~axiom as well. If not, $\a$ is not
principal and we use the induction hypothesis and permute contraction one step above or,
in the~case of $\Box_{GL}$ rule, we apply the~rule so that
the~conclusion is weakened by only one occurrence of $\a$.
For $\a$ not atomic and not boxed we consider its main connective and  use the~height preserving
invertibility of the~appropriate rule and by i.h. we apply
contraction on formula(s) of lower weight and then the~rule again.

Let us therefore only spell out the step for a formula of the form $\a=\Box\be$. If neither of its two occurrences  is principal in the last step of the proof of the premise of the contraction, we simply permute the weakening upwards and use the induction hypothesis. So consider one occurrence of $\Box\be$ is the principle formula of a $\Box_{GL}$ inference.
\begin{itemize}

\item[-] Contraction-right on $\Box\beta$, with one of the occurrences of $\Box\beta$ principle of a~$\Box_{GL}$ inference: we use the~$\Box_{GL}$ rule so that we do not weaken by the~other occurrence of $\Box\beta$ in the~conclusion.

\item[-] Contraction-left on $\Box\beta$, with the occurrences of $\Box\beta$ principle of a~$\Box_{GL}$ inference: we permute the~proof as follows and use a contraction on a simpler formula, h.p. admissible by the induction hypothesis:

\begin{prooftree} 
\AxiomC{$\be,\be,\Gamma,\Box\g\Rightarrow \g$}
\RightLabel{$\Box_{GL}\quad\Longrightarrow$}
\UnaryInfC{$\Box\be,\Box\be,\Box\Gamma,\Pi\Rightarrow\Box\g
,\Sigma$}
\RightLabel{c-l}
\UnaryInfC{$\Box\be,\Box\Gamma,\Pi\Rightarrow\Box\g ,\Sigma$}
\AxiomC{$\be,\be,\Gamma\Rightarrow \g$} 
\RightLabel{c-l, i.h.}
\UnaryInfC{$\be,\Gamma,\Box\g\Rightarrow\g$}
\RightLabel{$\Box_{GL}$}
\UnaryInfC{$\Box\be,\Box\Gamma,\Pi\Rightarrow\Box\g,\Sigma$}
\noLine\BinaryInfC{ }
\end{prooftree} 
\end{itemize}
All the above permutations are easily seen, using the induction hypothesis, to be height-preserving.
\end{proofof}

\subsection{Terminating proof-search in $G_{GL}$}\label{sub:termGL}

The proof-search strategy we adopt is based on applying the rules of $G_{GL}$ backwards to a given sequent, so that we always first apply the invertible rules and then, when it is no longer possible and if we haven't reached an axiom or a sequent with no boxed formulas on the right, we perform a modal jump --- we apply the $\Box_{GL}$ rule backwards. We prefer to pack all the invertible steps into a single step, therefore it is useful to define the following notions of a critical sequent and a closure of a sequent first:

\begin{defn}\label{def:cl}
A sequent is called \emph{critical}, if no invertible rule can be applied to it backwards. For a~sequent $(\Gamma\Rightarrow\Delta)$,
consider the~smallest set of sequents containing $(\Gamma\Rightarrow\Delta)$ and closed under backward
applications of the~invertible rules of $G_{GL}$. The~\emph{closure} of a~sequent $(\Gamma\Rightarrow\Delta)$, denoted $Cl(\Gamma;\Delta)$, is then the~subset of all \emph{critical} sequents contained in the set.
\end{defn}

Note that the closure of any sequent is finite, and that a critical sequent is of the form $(\Pi,\Box\Gamma\Rightarrow\Box\Delta,\Lambda)$, with $\Pi,\Lambda$ multisets of atomic formulas, and its closure is the singleton of the sequent itself. 

For a sequent $S= (\Gamma\Rightarrow\Delta)$ and finite multisets $\Theta,\Omega$, let $S(\Theta;\Omega)$ denote the sequent $\Gamma,\Theta\Rightarrow\Delta,\Omega$. 
The closure satisfies the following lemma, proof of which is immediate from the definition of the closure:

\begin{lem}\label{lem:cl}
Let $S$
be a~sequent, $Cl(S)=\{S_1,\ldots,S_n\}$, and $\Theta,\Omega$ arbitrary finite multisets of formulas.
Then:
\begin{itemize}
\item[-]
$S_1,\ldots,S_n\vdash_{G_{GL}}S$\\
         if $\ \vdash_{G_{GL}}S$ then
         $\ \vdash_{G_{GL}}S_i$ for each i.

\medskip
\item[-]        $S_1(\Theta;\Omega),\ldots,S_n(\Theta;\Omega)\vdash_{G_{GL}}S(\Theta;\Omega)$\\
         if $\ \vdash_{G_{GL}}S(\Theta;\Omega)$ then
         $\ \vdash_{G_{GL}}S_i(\Theta;\Omega)$ for each i.
\end{itemize}
\end{lem}

The proof-search procedure can now be described by creating a proof-search tree as follows: we start with creating a root and labeling it with the given sequent. For every node we have created, we proceed as follows: if it is labeled with a non-critical sequent, we compute its closure, and create a child-node for each sequent in the closure and label it with the sequent (thus creating a finite conjunctive branching). If a node is labeled with a critical sequent of the form $(\Pi,\Box\Gamma\Rightarrow\Box\Delta,\Lambda)$, we distinguish the following cases: if $\Pi\cap\Lambda\neq\emptyset$ or $\bot\in\Pi$ or $\Gamma\cap\Delta\neq\emptyset$ we mark the node a provable leaf, if it is not the case and $\Delta=\emptyset$ we mark the node an unprovable leaf, and in the remaining case we apply the $\Box_{GL}$ rule backwards: we create $|\Delta|$ children nodes and label them, for each $\Box\a\in\Delta$, with the premise of the $\Box_{GL}$ inference with $\Box\a$ principal (thus creating a finite disjunctive branching). 

Checking whether $\Gamma\cap\Delta\neq\emptyset$ before applying the modal rule backwards works as a~simple loop preventing mechanism --- we do not apply the~$\Box_{GL}$ rule backwards with $\Box\a$ principle if the~diagonal
formula $\Box\a$ is already in the~antecedent (in which case the sequent in question is clearly provable). This is crucial since it enables us to bound the~number of $\Box_{GL}$ inferences along each branch, and consequently also the weight of sequents occurring in a~proof search for a fixed sequent.

\begin{lem} \label{lem:GLterm} Proof search in the~calculus $G_{GL}$ always
terminates.
\end{lem}

\begin{proofof}{Lemma \ref{lem:GLterm}}

Consider a~proof search for a~sequent $(\Phi\Rightarrow\Psi)$. Let
$n$ be the~number of boxed subformulas contained in multisets $\Phi,\Psi$. This is, by the subformula property, by the nature of the~$\Box_{GL}$ rule, and by the loop-preventing mechanism described above, an upper bound of the 
number of $\Box_{GL}$ inferences along a~branch of the~proof search
tree. The reason is that, along a given branch, we never apply the~$\Box_{GL}$ rule with the same principal formula twice. In general, the weight of sequents increases whenever the~$\Box_{GL}$ rule is applied backwards. But the bound on the number of such steps along a single branch enables us to give an upper bound on the~weight of sequents occurring in the~fixed proof
search: namely, $c=2^{n}w(\Phi,\Psi)$ is an~upper bound of the~weight of a~sequent occurring in the~proof search for a~sequent
$(\Phi\Rightarrow\Psi)$. 

For any multiset $\Gamma$ occurring in the proof-search tree, let now $b(\Gamma)$ denote the~number of boxed formulas in $\Gamma$ counted
as a~set. For a~sequent $(\Gamma\Rightarrow\Delta)$ occurring during the proof-search, consider an~ordered pair $\langle c - b(\Gamma), w(\Gamma ,\Delta) \rangle$. This measure strictly
decreases in every backward application of a~rule in terms of
the~lexicographical ordering: $c$ is certainly greater or equal to the~maximal number of boxed formulas in the~antecedent which can occur during the~proof search, so the~first number does not decrease below zero. When an~invertible rule is applied backwards, the~weight of a~sequent strictly decreases, therefore for a non-critical sequent, all sequents from its closure are of strictly smaller weight. When the~$\Box_{GL}$ rule is applied backwards, $b$ increases, and so $c - b$ decreases, therefore the measure decreases. \footnote{Another way (closer to the~approach of \cite{Heu96} or \cite{Heu98}) how to formulate a measure is the~following: for a~sequent $(\Gamma;\Delta)$ consider the~function
$
f(\Gamma;\Delta)= c^{2}-cb(\Gamma)+w(\Gamma ,\Delta).
$
The function (values of which are non-negative integers) decreases
in every backward application of a~rule in a~proof search for
$(\Phi\Rightarrow\Psi)$. ($c^{2}$ is included to ensure that $f$ doesn't decrease below zero, and $cb(\Gamma)$ balances
the~possible increase of $w(\Gamma\Rightarrow\Delta)$ in the~case
of a~backward application of the~$\Box_{GL}$-rule.)}

\end{proofof}

\subsubsection{Extracting a proof} From a proof-search tree for a provable sequent we are expected to be able to extract an actual proof of the sequent. The tree is finite, and all the leaves are marked either provable, or not provable. We can extend the marking in an obvious way to all the nodes: if a node is a parent node of a conjunctive branching, we mark it provable if and only if all its children are marked provable, if a node is a parent node of a disjunctive branching, we mark it provable if and only if at least one of its children is marked provable. If the root is marked provable, the sequent we started with has a proof presented by a tree of nodes marked provable, generated by the root. It is a routine induction to see it is indeed a proof.

\subsection{Sequent calculi for \textbf{Grz}}
We present a calculus for Grzegorczyk's logic with a loop-preventing mechanism built into the syntax of sequents. Namely we include a third multiset of boxed formulas in a sequent, thus sequents are now of the form $\Box\Sigma|\Gamma\Rightarrow\Delta$. The third multiset is used to store the boxed formulas of the $\Box_{T}^{+}$ inferences and the diagonal formulas of the $\Box_{Grz2}^{+}$ inferences when the rules are applied backwards to prevent unnecessary looping. This strategy was inspired by work of Heuerding \cite{Heu96,Heu98}.

To improve readability, we denote the diagonal formula $\Box(\a\rightarrow\Box \a)$ of the Grzegorczyk's axiom by $D(\a)$ in the following text.

\begin{defn} Sequent calculus $G_{Grz}^{+}$:

\begin{prooftree}
\AxiomC{$\Box\Sigma |\Gamma,p\Rightarrow p ,\Delta$}
\AxiomC{$\Box\Sigma |\Gamma,\bot\Rightarrow\Delta$}
\noLine\BinaryInfC{ }
\end{prooftree}

\begin{prooftree}
\AxiomC{$\Box\Sigma |\Gamma,\a,\be\Rightarrow\Delta$}
\RightLabel{$\wedge$-l}
\UnaryInfC{$\Box\Sigma |\Gamma,\a\wedge\be\Rightarrow\Delta$}
\AxiomC{$\Box\Sigma |\Gamma\Rightarrow\a,\be,\Delta$}
\RightLabel{$\vee $-r}
\UnaryInfC{$\Box\Sigma|\Gamma\Rightarrow\a\vee\be,\Delta$}
\noLine\BinaryInfC{ }
\end{prooftree}

\begin{prooftree}
\AxiomC{$\Box\Sigma |\Gamma,\a\Rightarrow\Delta$}
\RightLabel{$\neg$-r}
\UnaryInfC{$\Box\Sigma |\Gamma\Rightarrow\neg\a,\Delta$}
\AxiomC{$\Box\Sigma |\Gamma\Rightarrow\a,\Delta$}
\RightLabel{$\neg$-l}
\UnaryInfC{$\Box\Sigma |\Gamma,\neg\a\Rightarrow\Delta$}
\noLine\BinaryInfC{  }
\end{prooftree}

\begin{prooftree} 
\AxiomC{$\Box\Sigma |\Gamma\Rightarrow\a,\Delta$}
\AxiomC{$\Box\Sigma |\Gamma\Rightarrow\be,\Delta$}
\RightLabel{$\wedge$-r}
\BinaryInfC{$\Box\Sigma |\Gamma\Rightarrow\a\wedge\be,\Delta$}
\AxiomC{$\Box\Sigma |\Gamma,\a\Rightarrow\Delta$}
\AxiomC{$\Box\Sigma |\Gamma,\be\Rightarrow\Delta$}
\RightLabel{$\vee$-l}
\BinaryInfC{$\Box\Sigma | \Gamma ,\a\vee\be\Rightarrow\Delta$}
\noLine\BinaryInfC{  }
\end{prooftree}

\begin{prooftree} 
\AxiomC{$\Box\a,\Box\Sigma |\Gamma,\a\Rightarrow\Delta$}
\RightLabel{$\Box_{T}^{+}$}
\UnaryInfC{$\Box\Sigma|\Gamma,\Box\a\Rightarrow\Delta$}
\end{prooftree} 

\begin{prooftree}
\AxiomC{$\Box\Gamma|\emptyset\Rightarrow\a$}
\RightLabel{$\Box_{Grz1}^{+},D(\a)\in\Box\Gamma$}
\UnaryInfC{$\Box\Gamma|\Pi\Rightarrow\Box\a,\Delta$}
\AxiomC{$\Box\Gamma,D(\a)|\Gamma\Rightarrow\a$}
\RightLabel{$\Box_{Grz2}^{+},D(\a)\notin\Box\Gamma$}
\UnaryInfC{$\Box\Gamma|\Pi\Rightarrow\Box\a,\Delta$}
\noLine\BinaryInfC{}
\end{prooftree}
\end{defn}

In the~$\Box_{Grz}^{+}$ rules, $\Pi$ contains only propositional
variables and $\Delta$ contains only propositional variables and
boxed formulas. The notion of a principal formula is similar to the previous case. All the propositional rules are easily seen to be height-preserving invertible, there is one additional invertible rule here:

\begin{lem}\label{lem:Tinv}
$\Box_{T}^{+}$ is height-preserving invertible.
\end{lem}

\begin{proofof}{Lemma \ref{lem:Tinv}}
Proof is a routine induction on the height of the proof of the premise and we leave it to the reader.
\end{proofof}

\begin{lem}\label{lem:wGrz+} Weakening rules are
admissible in $G_{Grz}^{+}$.
\end{lem}
\begin{proofof}{Lemma \ref{lem:wGrz+}} 
The~weakening rules we consider in this paper are: 

\begin{prooftree} 
\AxiomC{$\Box\Sigma|\Gamma\Rightarrow \Delta$} 
\RightLabel{w-l}
\UnaryInfC{$\Box\Sigma|\Gamma,\a\Rightarrow\Delta$}
\AxiomC{$\Box\Sigma|\Gamma \Rightarrow \Delta$}
\RightLabel{w-r}
\UnaryInfC{$\Box\Sigma|\Gamma \Rightarrow\Delta,\a$}
\AxiomC{$\Box\Sigma|\Gamma \Rightarrow \Delta$}
\RightLabel{w-l+}
\UnaryInfC{$\Box\Sigma, \Box \a|\Gamma\Rightarrow\Delta$}
\noLine\TrinaryInfC{ }
\end{prooftree}
The~proof is by induction on the~weight of the~principal weakening formula $\a$ and, for each weight, on the~height of the~proof of the~premise. We prove admissibility of the three weakening rules simultaneously. 

For an~atomic weakening formula the~proof is obvious - note that
atomic weakening is built in axioms as well as in the~$\Box_{Grz}^{+}$-rules. (The case of weakening-r by $\Box\a$ when the last inference is one of the~$\Box_{Grz}^{+}$-rules is then also obvious since it is built-in the~rules as well.)
For non atomic and not boxed formula we consider its main connective and use height-preserving
invertibility of the~corresponding propositional rule, weaken by formula(s) of
lower weight (admissible by the induction hypothesis), and then apply the~rule. We next consider the~weakening formula being of the~form $\Box\be$. 

\medskip\par\noindent
\textbf{Weakening-right:} consider the last step of the proof of the premise of the weakening inference. If it is a $\Box_{Grz}^+$ rule, we can use the rule so that the weakening by $\Box\be$ is built-in its conclusion. If it is an invertible rule, we permute the weakening upwards.

\medskip\par\noindent
\textbf{Weakening-left:} consider the last step of the proof of the premise of the weakening inference. If it is an invertible rule, we permute the weakening upwards. Let us consider last step is a $\Box_{Grz}^+$ rule, then we permute as follows:
\begin{prooftree}
\AxiomC{$D(\g),\Box\Sigma|\emptyset\Rightarrow\g$}
\RightLabel{$\Box_{Grz1}^{+}\ \Longrightarrow$}
\UnaryInfC{$D(\g),\Box\Sigma|\Pi\Rightarrow\Box\g,\Delta$}
\RightLabel{w-l}
\UnaryInfC{$D(\g),\Box\Sigma|\Box\be,\Pi\Rightarrow\Box\g,\Delta$}
\AxiomC{$D(\g),\Box\Sigma|\emptyset\Rightarrow\g$}
\RightLabel{w-l$^+$}
\UIC{$D(\g),\Box\Sigma,\Box\be|\emptyset\Rightarrow\g$}
\RightLabel{$\Box_{Grz1}^{+}$}
\UnaryInfC{$D(\g),\Box\Sigma,\Box\be|\Pi\Rightarrow\Box\g,\Delta$}
\RightLabel{w-l} 
\UnaryInfC{$D(\g),\Box\Sigma,\Box\be|\be,\Pi\Rightarrow\Box\g,\Delta$}
\RightLabel{$\Box_T^{+}$}
\UnaryInfC{$D(\g),\Box\Sigma|\Box\be,\Pi\Rightarrow\Box\g,\Delta$}
\noLine\BIC{ }
\end{prooftree}

\begin{prooftree}
\AxiomC{$D(\be),\Box\Sigma|\Sigma\Rightarrow\be$}
\RightLabel{$\Box_{Grz2}^{+}\ \Longrightarrow$}
\UnaryInfC{$\Box\Sigma|\Pi\Rightarrow\Box\be,\Delta$}
\RightLabel{w-l}
\UnaryInfC{$\Box\Sigma|\Box\a,\Pi\Rightarrow\Box\be,\Delta$}
\AxiomC{$D(\be),\Box\Sigma|\Sigma\Rightarrow\be$}
\RightLabel{w-l$^+$,l}
\UIC{$D(\be),\Box\Sigma,\Box\a|\a,\Sigma\Rightarrow\be$}
\RightLabel{$\Box_{Grz2}^{+}$}
\UnaryInfC{$\Box\Sigma,\Box\a|\Pi\Rightarrow\Box\be,\Delta$}
\RL{w-l}
\UIC{$\Box\Sigma,\Box\a|\a,\Pi\Rightarrow\Box\be,\Delta$}
\RightLabel{$\Box_T^{+}$}
\UnaryInfC{$\Box\Sigma|\Box\a,\Pi\Rightarrow\Box\be,\Delta$}
\noLine\BIC{ }
\end{prooftree}

\par\noindent
\textit{Remark:} The two transformations above are clearly not height-preserving, therefore weakenings are in general not height-preserving admissible. However, one can show, that weakening rules with $\a$ principal are admissible and the height only increases by  the box depth $d(\a)$.

\medskip\par\noindent
\textbf{Weakening-l+:} Notice that w-l+ is built in the axioms. If the~last inference of the proof of the premise of the weakening is a~propositional inference or a~$\Box_T^{+}$ inference, we just use the~i.h., a weakening one step above, and use the~appropriate rule again.

Let the~last inference of the proof of the premise of the weakening be
a~$\Box_{Grz1}^{+}$ inference, w-l+ permutes over the inference as follows:

\begin{prooftree}
\AxiomC{$D(\g),\Box\Sigma|\emptyset\Rightarrow\g$}
\RightLabel{$\Box_{Grz1}^{+}\ \Longrightarrow$}
\UnaryInfC{$D(\g),\Box\Sigma|\Pi\Rightarrow\Box\g,\Delta$}
\RightLabel{w-l$^{+}$}
\UnaryInfC{$D(\g),\Box\Sigma,\Box\be|\Pi\Rightarrow\Box\g,\Delta$}
\AxiomC{$D(\g),\Box\Sigma|\emptyset\Rightarrow\g$}
\RightLabel{w-l$^{+}$}
\UnaryInfC{$D(\g),\Box\Sigma,\Box\be|\emptyset\Rightarrow\g$}
\RightLabel{$\Box_{Grz1}^{+}$}
\UnaryInfC{$D(\g),\Box\Sigma,\Box\be|\Pi\Rightarrow\Box\g,\Delta$}
\noLine\BIC{ }
\end{prooftree}
Let the~last inference be
a~$\Box_{Grz2}^{+}$ inference, w-l+ permutes over the inference as follows:

\begin{prooftree}
\AxiomC{$D(\g),\Box\Sigma|\Sigma\Rightarrow\g$}
\RightLabel{$\Box_{Grz2}^{+}\ \Longrightarrow$}
\UnaryInfC{$\Box\Sigma|\Pi\Rightarrow\Box\g,\Delta$}
\RightLabel{w-l$^{+}$}
\UnaryInfC{$\Box\Sigma,\Box\be|\Pi\Rightarrow\Box\g,\Delta$}
\AxiomC{$D(\g),\Box\Sigma|\Sigma\Rightarrow\g$}
\RightLabel{w-l,l$^{+}$}
\UnaryInfC{$D(\g),\Box\Sigma,\Box\be|\be,\Sigma\Rightarrow\g$}
\RightLabel{$\Box_{Grz2}^{+}$}
\UnaryInfC{$\Box\Sigma,\Box\be|\Pi\Rightarrow\Box\g,\Delta$}
\noLine\BIC{ }
\end{prooftree}
The two permutations above are in fact height-preserving. 
\end{proofof}

\begin{lem}\label{lem:contrGrz+} Contraction rules are
height-preserving admissible in $G_{Grz}^{+}$.
\end{lem}
\begin{proofof}{Lemma \ref{lem:contrGrz+}} 

The~contraction rules are:
\begin{prooftree}
\AxiomC{$\Sigma|\Gamma,\a,\a\Rightarrow \Delta$}
\RightLabel{c-l}
\UnaryInfC{$\Sigma|\Gamma,\a\Rightarrow\Delta$}
\AxiomC{$\Sigma|\Gamma \Rightarrow \Delta,\a,\a$} 
\RightLabel{c-r}
\UnaryInfC{$\Sigma|\Gamma\Rightarrow\Delta,\a$}
\AxiomC{$\Sigma,\Box\a,\Box \a|\Gamma\Rightarrow \Delta$}
\RightLabel{c-l+}
\UnaryInfC{$\Sigma,\Box\a|\Gamma\Rightarrow\Delta$}
\noLine\TIC{ }
\end{prooftree}

The~proof is by induction on the~weight of the~contraction formula
and, for each weight, on the~height of the~proof of the~premise.
The~induction runs simultaneously for all the~contraction rules.
We use the~height preserving invertibility of rules. Note that in
the~contr-l+ rule the~contraction formula is always of the~form
$\Box\a$.

For $\a$ atomic, if the~premise is an~axiom,
the~conclusion is an~axiom as well. If not, $\a$ is not
principal and we use the induction hypothesis and apply contraction one step above or,
in the~case of $\Box_{Grz}^{+}$ rules, we apply the~rule so that
the~conclusion is weakened by only one occurrence of $\a$.
For $\a$ not atomic and not boxed we consider its main connective and use the~height preserving
invertibility of the~corresponding propositional rule and by the induction hypothesis we apply
contraction on subformula(s) of lower weight and then the~rule again.
The~third multiset does not make any difference here and it works
precisely as in the~classical logic. All the~steps described so far are height preserving.

Now suppose the~contraction formula to be of the~form $\Box\be$. We
distinguish the following cases:

\medskip\par\noindent
(i) Both occurrences of the~contraction formula are principal of
a~$\Box_{Grz1}^{+}$ inference in the~antecedent, in this case the only possibility is c-l$^+$. Then we permute the
proof as follows using the~the induction hypothesis:

\begin{prooftree} 
\AxiomC{$\Box\Sigma,\Box\be,\Box\be|\emptyset\Rightarrow\g$}
\RightLabel{$\Box_{Grz1}^{+}\quad\Longrightarrow$}
\UnaryInfC{$\Box\Sigma,\Box\be,\Box\be|\Pi\Rightarrow\Box\g
,\Delta$}
\RightLabel{c-l+}
\UnaryInfC{$\Box\Sigma,\Box\be|\Pi\Rightarrow\Box\g,\Delta$}
\AxiomC{$\Box\Sigma,\Box\be,\Box\be|\emptyset\Rightarrow\g$}
\RightLabel{c-l+}
\UnaryInfC{$\Box\Sigma,\Box\be|\emptyset\Rightarrow\g$}
\RightLabel{$\Box_{Grz1}$}\UnaryInfC{$\Box\Sigma,\Box\be|\Pi\Rightarrow\Box\g,\Delta$}
\noLine\BinaryInfC{ }
\end{prooftree} 
In this case, $D(\g)\in\Box\Sigma$, or $D(\g)=\Box\be$. The permutation is obviously height preserving.

\medskip\par\noindent
(ii) Both occurrences of the~contraction formula are principal of
a~$\Box_{Grz2}^{+}$ inference in the~antecedent, again c-l$^+$ is the only possibility. Then we permute the proof as follows:

\begin{prooftree} 
\AxiomC{$D(\g),\Box\Sigma,\Box\be,\Box\be|\Sigma,\be,\be\Rightarrow\g$}
\RightLabel{$\Box_{Grz2}^{+}\quad\Longrightarrow$}
\UnaryInfC{$\Box\Sigma,\Box\be,\Box\be|\Pi\Rightarrow\Box\g
,\Delta$}
\RightLabel{c-l+}
\UnaryInfC{$\Box\Sigma,\Box\be|\Pi\Rightarrow\Box\g,\Sigma$}
\AxiomC{$D(\g),\Box\Sigma,\Box\be,\Box\be|\Sigma,\be,\be\Rightarrow\g$}
\RightLabel{c-l,l+}
\UnaryInfC{$D(\g),\Box\Sigma,\Box\be|\Sigma,\be\Rightarrow\g$}
\RightLabel{$\Box_{Grz2}$}
\UnaryInfC{$\Box\Sigma,\Box\be|\Pi\Rightarrow\Box\g,\Sigma$}
\noLine\BinaryInfC{ }
\end{prooftree} 
Here, $D(\g)\notin\Box\Sigma$ and $D(\g)\neq\Box\be$. The permutation is obviously height preserving. 

\medskip\par\noindent
(iii) One occurrence of the~contraction formula is the~principal formula of
a~$\Box_{T}^{+}$ inference in the~antecedent. Then we permute the
proof as follows using the~i.h.\ and the~height preserving
invertibility of the~$\Box_{T}^{+}$ rule:

\begin{prooftree}
\AxiomC{$\Sigma,\Box\be|\Box\be,\be,\Gamma\Rightarrow\Delta$}
\RightLabel{$\Box_{T}^{+}\quad\Longrightarrow$}
\UnaryInfC{$\Sigma|\Box\be,\Box\be,\Gamma\Rightarrow\Delta$}
\RightLabel{c-l}
\UnaryInfC{$\Sigma|\Box\be,\Gamma\Rightarrow\Delta$}
\AxiomC{$\Sigma,\Box\be|\Box\be,\be,\Gamma\Rightarrow\Delta$}
\RightLabel{invert.}
\UnaryInfC{$\Sigma,\Box\be,\Box\be|\be,\be,\Gamma\Rightarrow\Delta$}
\RightLabel{c-l,l+}
\UnaryInfC{$\Sigma,\Box\be|\be,\Gamma\Rightarrow\Delta$}
\RightLabel{$\Box_{T}^{+}$}
\UnaryInfC{$\Sigma|\Box\be,\Gamma\Rightarrow\Delta$}
\noLine\BinaryInfC{ }
\end{prooftree} 
The permutation is height preserving since the~steps c-l, c-l+, and invert. do not increase the~height of the~proof.

\medskip\par\noindent
(iv) One occurrence of the~contraction formula is the~principal formula in
the~succedent and the last inference is a~$\Box_{Grz}^{+}$ inference. 
Then we use the~$\Box_{Grz}^{+}$ rule so
that the~conclusion is not weakened by the~other occurrence of $\Box\be$. This step is obviously height preserving.
If the~contraction formula is not the~principal formula and
the~last step is a~$\Box_{Grz}^{+}$ inference, $\Box\be$ is in
$\Delta$. Then we use the~$\Box_{Grz}^{+}$ rule so that
the~conclusion is weakened by only one occurrence of
the~contraction formula. If the~last step is another inference, we
use contraction one step above on the~proof of lower height. If it
is an~axiom, the~conclusion of the~desired contraction
is an~axiom as well. Again, all the~steps are height
preserving.
\end{proofof}

Next we want to relate the calculus $G_{Grz}^{+}$ to the standard sequent calculus $G_{Grz}$ of \cite{Avr84}. For this, we consider a multiset variant of the latter. The calculus is known to be complete, and the rules of weakening are easily proved to be admissible by a similar argument that used in Lemma \ref{lem:strGL}.

\begin{defn}
Sequent calculus $G_{Grz}$ results from the non-modal part of the calculus $G_{GL}$ adding the~following two modal rules:

\begin{prooftree}
\AxiomC{$\Gamma,\Box\a,\a\Rightarrow \Delta$}
\RightLabel{$\Box_{T}$}
\UnaryInfC{$\Gamma,\Box\a\Rightarrow\Delta$}
\AxiomC{$\Box\Gamma,D(\a)\Rightarrow\a$}
\RightLabel{$\Box_{Grz}$}
\UnaryInfC{$\Box\Gamma,\Pi\Rightarrow\Box\a,\Delta$}
\noLine\BIC{ }
\end{prooftree}
\end{defn}

\begin{lem}\label{lem:eqGrz+} The~calculi $G_{Grz}$ and
$G_{Grz}^{+}$ are equivalent:
$$\vdash_{G_{Grz}}\Gamma\Rightarrow\Delta\quad \mbox{iff}
\quad\vdash_{G_{Grz}^{+}}\emptyset|\Gamma\Rightarrow\Delta. $$
\end{lem}

\begin{proofof}{Lemma \ref{lem:eqGrz+}}
$\mbox{}$

\par\noindent
\textit{The~right-left implication:} deleting the~"$|$" symbol from
a~$G_{Grz}^{+}$ proof of $(\emptyset|\Gamma\Rightarrow\Delta)$
yields correct instances of rules of $G_{Grz}$, except
the~$\Box_{Grz1}^{+}$ rule. It has to be simulated as follows:

\begin{prooftree}
\AxiomC{$\Box\Gamma,D(\a)\Rightarrow\a$}
\RightLabel{$\Box_{Grz}$}
\UnaryInfC{$\Box\Gamma,\Pi\Rightarrow\Box\a,\Delta$}
\RightLabel{admiss. w-l}
\UnaryInfC{$\Box\Gamma,D(\a),\Pi\Rightarrow\Box\a,\Delta$}
\end{prooftree}
We end up with a~$G_{Grz}$ proof of $\Gamma\Rightarrow\Delta$. (Lemma \ref{lem:cutG} below states that the calculus $G_{Grz}$ is complete, and soundness of weakening entails that weakening is indeed admissible in $G_{Grz}$.) 

\medskip\par\noindent
\textit{The~left-right implication:} the idea is to add a third, empty multiset to all the sequents in a proof. This yields correct instances of the axioms as well as the~propositional rules. The~$\Box_{T}$ 
rule has to be simulated as follows, using invertibility of $\Box_{T}^{+}$ rule and admissibility of contraction: 

\begin{prooftree}
\AxiomC{$\emptyset|\Box\a,\a\Gamma\Rightarrow\Delta$}
\RL{inv. of $\Box_{T}^{+}$}
\UIC{$\Box\a|\a,\a,\Gamma\Rightarrow\Delta$}
\RL{c-l}
\UIC{$\Box\a|\a,\Gamma\Rightarrow\Delta$}
\RL{$\Box_{T}^{+}$}
\UIC{$\emptyset|\Box\a,\Gamma\Rightarrow\Delta$}
\end{prooftree}

\par\noindent
The~$\Box_{Grz}$ rule is simulated as follows ($D(\a)\notin\Box\Gamma$):

\begin{prooftree}
\AxiomC{$\emptyset|\Box\Gamma,\Box(\neg\a\vee\Box\a)\Rightarrow\a$}
\RightLabel{inv. of $\Box_{T}^{+}$}
\UnaryInfC{$\Box\Gamma,\Box(\neg\a\vee\Box\a)|\Gamma,\neg\a\vee\Box\a\Rightarrow\a$}
\RightLabel{inv. of $\vee$-l and $\neg$-l}
\UnaryInfC{$\Box\Gamma,\Box(\neg\a\vee\Box\a)|\Gamma\Rightarrow\a,\a$}\RightLabel{admiss. c-r}
\UnaryInfC{$\Box\Gamma,\Box(\neg\a\vee\Box\a)|\Gamma\Rightarrow\a$}
\RightLabel{$\Box_{Grz2}^{+}$}
\UnaryInfC{$\Box\Gamma|\Pi\Rightarrow\Box\a,\Delta$}
\RightLabel{admiss. w-l inferences}
\UnaryInfC{$\Box\Gamma|\Gamma,\Pi\Rightarrow\Box\a,\Delta$}
\RightLabel{$\Box_{T}^{+}$ inferences}
\UnaryInfC{$\emptyset|\Box\Gamma,\Pi\Rightarrow\Box\a,\Delta $}
\end{prooftree}
If $D(\a)\in\Gamma $, we use some admissible
c-l+ inferences before the~$\Box_{Grz1}^{+}$ inference is
used.
\end{proofof}

\subsection{Terminating proof-search in $G_{Grz}^{+}$}
We will restrict ourselves to proof-search for sequents of the form $(\emptyset|\Phi\Rightarrow\Psi)$ with the third multiset empty and $\Phi$ and $\Psi$ arbitrary finite multisets of formulas (Lemma \ref{lem:eqGrz+} justifies this restriction). 

The notion of critical sequent and the closure of a sequent for $G_{Grz}^{+}$ is the same as given in Definition \ref{def:cl}, only with a third multiset added. Recall that now also the $\Box_T^+$ rule is invertible, so critical sequents are of the form:
$
\Box\Gamma|\Pi\Rightarrow\Box\Delta,\Lambda
$
with $\Pi,\Lambda$ atomic.

For a sequent $S=(\Box\Sigma|\Gamma\Rightarrow\Delta)$ and finite multisets $\Lambda,\Theta,\Omega$, let $S(\Lambda|\Theta;\Omega)$ denote the sequent $(\Box\Sigma,\Box\Lambda|\Gamma,\Theta\Rightarrow\Delta,\Omega)$. 
The closure satisfies the following lemma, essentially the same as Lemma \ref{lem:cl}:

\begin{lem}\label{lem:cl+}
Let $S$
be a~sequent, $Cl(S)=\{S_1,\ldots,S_n\}$, and $\Lambda|\Theta,\Omega$ arbitrary finite multisets of formulas.
Then:
\begin{itemize}
\item[-]
$S_1,\ldots,S_n\vdash_{G_{GL}}S$\\
         if $\ \vdash_{G_{GL}}S$ then
         $\ \vdash_{G_{GL}}S_i$ for each i.

\medskip
\item[-]        $S_1(\Lambda|\Theta;\Omega),\ldots,S_n(\Lambda|\Theta;\Omega)\vdash_{G_{GL}}S(\Lambda|\Theta;\Omega)$\\
         if $\ \vdash_{G_{GL}}S(\Lambda|\Theta;\Omega)\ $ then
         $\ \vdash_{G_{GL}}S_i(\Lambda|\Theta;\Omega)$ for each i.
\end{itemize}
\end{lem}

\medskip\par\noindent
Before we continue to describe  proof-search and prove its termination, we briefly discuss forms of looping we prevent by using the specific form of the calculus.
\medskip
\subsubsection{Reflexive looping}
This simple looping occurs when one searches for proofs in the calculus $G_{Grz}$ and applies the $\Box_T$ rule backwards repeatedly with the same principal formula. Such looping is prevented by the presence of the third storage multiset in sequents and by the particular form of $\Box_T^+$ rule we use --- when this rule is applied backwards, it remembers that the principle formula has already been treated.

\subsubsection{Transitive looping}
Another looping phenomenon arises when one tries to search for a proof of the sequent $\Box\neg\Box
p\Rightarrow\Box p$ in the calculus $G_{Grz}$ --- it loops on the sequent $$
\Box\neg\Box
p,D(p)\Rightarrow p,\Box p.
$$
Such looping can be avoided and the diagonal formula plays a~crucial role here as a~natural loop-preventing mechanism again. We have made this mechanism explicit by splitting the~$\Box_{Grz}$ rule into two
cases distinguishing if the~diagonal formula is present in
the~antecedent or not. Consider the~$\Box_{Grz1}^{+}$ rule bottom up. When
the~diagonal formula is already in the~third multiset, we apply
the~rule so that we neither add the~diagonal formula to the~third
multiset, nor we add $\Gamma$ to the~antecedent.

\medskip\par\noindent
The proof-search procedure for Grzegorczyk logic is fully analogous to that for logic \textbf{GL}: we create a proof-search tree, using the strategy of alternating the closure step for non-critical sequents and a modal jump step for critical sequents. Also the labeling and extraction of an actual proof is carried out similarly.

\begin{lem}\label{lem:Grzterm} Proof search in $Gm_{Grz}^{+}$ for sequents of the~form
$(\emptyset|\Phi\Rightarrow\Psi)$ always terminates.
\end{lem}
\begin{proofof}{Lemma \ref{lem:Grzterm}}
Consider a~proof search for a~sequent
$(\emptyset|\Phi\Rightarrow\Psi)$. Let $n$ be the number of boxed subformulas occurring in the sequent
$(\emptyset|\Phi\Rightarrow\Psi)$. This number, as in the~case of
\textbf{GL}, is an upper bound on the~number of the~$\Box_{Grz}^{+}$ rules applied backwards along one
branch of the~proof search tree. Each backward application of the~$\Box_{Grz2}^{+}$ rule adds a new boxed formula in the storage multiset, but also a $\Box_{T}^{+}$ rule does so during the closure steps. Therefore $n^2$ is an~upper bound of the~number of formulas stored in $\Sigma$ if we do not duplicate them and count them as a set. (If we allowed duplicate formulas in $\Sigma$, we would need an~exponential function of $n$.)

With each sequent $(\Box\Sigma|\Gamma\Rightarrow\Delta)$ occurring
during the~proof search, we associate an~ordered pair $\langle n^2-
|\Sigma^{\circ}|,w(\Gamma,\Delta) \rangle$. Therefore
the~first number does not decrease below zero. The measure obviously
decreases in every backward application of a~rule of the~calculus.
For the~$\Box_{Grz2}^{+}$ rule, $|\Sigma^{\circ}|$ increases and so
$n^2-|\Sigma^{\circ}|$ decreases, while for other rules the~weight
$w(\Gamma,\Delta)$  decreases.

\end{proofof}

\subsection{Cut admissibility via completeness}\label{ss:cutG}

We do not give a constructive proof of completeness of the two calculi without the cut rule in this paper. Such a proof can be established using a proof-search method described in the previous subsections. One can argue that, for any given sequent, the proof-search tree either yields a proof, or can be used to construct a (finite) counterexample. Instead, we state the~completeness without the cut rule, and refer for a proof to Avron \cite{Avr84} who proved that the~calculi $G_{GL}$ and $G_{Grz}$ are complete without the~cut rule w.r.t. their respective Kripke semantics. Then an~easy semantical argument of soundness of the~cut rule entails its admissibility. Lemma \ref{lem:eqGrz+} then yields admissibility in $G_{Grz}^{+}$ of the~cut rule we will use later in the~proof of Theorem \ref{thm:Grz}.

Proofs of completeness can be found in\cite{Avr84} for \textbf{GL} and Grzegorczyk's logic, and \cite{Sve00} or \cite{SV82} for \textbf{GL}, where redundancy of the~cut rule is established through a~decision procedure which either creates a~cut-free proof, or a~Kripke counterexample to a~given sequent. Although both authors use a~formulation via sets of formulas, observe, that a~cut-free proof with sets can be equivalently formulated using multisets and contraction rules, which are, as we have proved, admissible in our cut-free calculi. Equivalently, if a~sequent does not have a~cut-free proof in the~system based on multisets, its set-based counterpart sequent does not have a~cut-free proof in the~system based on sets.

\begin{lem}\label{lem:cutG}(Avron \cite{Avr84}:)
There are a~canonical Kripke model $(W,<)$  and a~canonical
valuation $V$ such that:

\begin{itemize}
    \item[-] $<$ is irreflexive and transitive
    \item[-]  for every $w\in W$, the~set $\{v|v<w\}$ is finite
    \item[-]  if $(\Gamma\Rightarrow\Delta)$ has no cut-free proof in $G_{GL}$, then there is a~$w\in W$ such that $w\Vdash_{V}\a$ for every $\a\in \Gamma$ and $w\nVdash_{V}\be$ for every $\be\in\Delta$.
\end{itemize}
There are a~canonical Kripke model $(W,\leq)$ and a~canonical
valuation $V$ such that:
\begin{itemize}
   \item[-]  $\leq$ partially orders W
   \item[-]  for every $w\in W$, the~set $\{v|v\leq w\}$ is finite
   \item[-]  if $(\Gamma\Rightarrow\Delta)$ has no cut-free proof in $G_{Grz}$, then there is a~$w\in W$ such that $w\Vdash_{V}\a$ for every $\a\in \Gamma$ and $w\nVdash_{V}\be$ for every $\be\in\Delta$.
\end{itemize}
\end{lem}

\begin{proofof}{Lemma \ref{lem:cutG}}
See \cite{Avr84}. The~canonical model is built from all saturated
sequents (sequents closed under subformulas) that have no cut-free proof in appropriate calculi. The~lemma entails completeness of $G_{GL}$ w.r.t. transitive, conversely well-founded Kripke models; and completeness of $G_{Grz}$ w.r.t. transitive, reflexive and conversely well-founded Kripke models.

\end{proofof}

\begin{cor}\label{cor:cutG} The~cut rule
\begin{prooftree}\AxiomC{$\Gamma
\Rightarrow \Delta,\gamma$}\AxiomC{$\gamma,\Pi\Rightarrow\Lambda$}
\BinaryInfC{$\Gamma,\Pi\Rightarrow\Delta,\Lambda$}
\end{prooftree}
 is admissible in $G_{GL}$ and $G_{Grz}$.
\end{cor}
\begin{proofof}{Cor. \ref{cor:cutG}}
It is easy to give a~semantic argument of soundness of the~cut rule.
Given a~counterexample of the~conclusion
$(\Gamma,\Pi\Rightarrow\Delta,\Lambda)$ of a~cut inference, there is
a~counterexample to one of its premises: consider
the~counterexample (W,R) and a~world $w\in W$ in it such that
$w\Vdash_{V} A$ for every $\a\in \Gamma\cup\Pi$ and $w\nVdash_{V}\be$
for some $\be\in\Delta\cup\Lambda$. For any formula $\g$ it is either
the~case that $w\Vdash_{V}\g$, and then $w$ refutes $(\g,\Pi
\Rightarrow \Lambda)$, or $w\nVdash_{V}\g$, and then $w$ refutes
$(\Gamma\Rightarrow\Delta,\g)$.

Now Lemma \ref{lem:cutG} (completeness of $G_{GL}$ and $G_{Grz}$)
entails admissibility of the~cut rule in the~calculi.
\end{proofof}

\begin{cor}\label{cor:cutGrz+} The~cut rule
\begin{prooftree}\AxiomC{$\emptyset|\Gamma
\Rightarrow \Delta,\gamma$}\AxiomC{$\emptyset|\gamma,\Pi\Rightarrow\Lambda$}
\BinaryInfC{$\emptyset|\Gamma,\Pi\Rightarrow\Delta,\Lambda$}
\end{prooftree}
 is admissible in  $G_{Grz}^{+}$.
\end{cor}
\begin{proofof}{Cor. \ref{cor:cutGrz+}}
Follows from Corollary \ref{cor:cutG} and Lemma \ref{lem:eqGrz+}.
\end{proofof}

\begin{rem}\label{rem:cut}
The~above cut rule cannot be replaced by the~expected form of cut:
\begin{prooftree}\AxiomC{$\Box\Sigma|\Gamma\Rightarrow \Delta ,\a$}
\AxiomC{$\Box\Theta|\a,\Pi \Rightarrow \Lambda$}
\RightLabel{cut',}
\BinaryInfC{$\Box\Sigma,\Box\Theta|\Gamma,\Pi\Rightarrow
\Delta,\Lambda$}
\end{prooftree}
since it is not admissible. The~counterexample is
the~following instance of cut':
\begin{prooftree}
\AxiomC{$\Box p,D(p)|p\Rightarrow
p$}
\RightLabel{$\Box_{Grz2}^{+}$}
\UnaryInfC{$\Box
p|\emptyset\Rightarrow\Box p$}
\AxiomC{$\Box p|p\Rightarrow
p$}
\RL{$\Box_T^+$}
\UIC{$\emptyset|\Box p\Rightarrow
p$}
\RightLabel{cut'}
\BinaryInfC{$\Box p|\emptyset\Rightarrow p$}
\end{prooftree}
which results in the sequent $(\Box p|\emptyset\Rightarrow p)$ unprovable
in $G_{Grz}^{+}$.
\end{rem}

\medskip\par\noindent
\textit{A fixed-point trick: }Before we proceed to the proof of uniform interpolation, we state a lemma which later will play a crucial role in a termination argument for the definition of the interpolants. The lemma is a simple example of a particular fixed point existence in modal logic \textbf{K4}. Namely, a recursive equivalence ${\diamondsuit((\a\vee x)\wedge\be)\equiv x}$ has a solution $x=\diamondsuit(\a\wedge\be)$. We will however need a bit more complicated form of the statement, which is the following lemma:
\begin{lem}\label{lem:k4trick}
Let $\de = \bigwedge\limits_{i}\a_{i}\wedge\be$, $i=1\ldots k$. Then the following sequents are provable:
$$
\vdash_{G_{GL}}\diamondsuit(\bigwedge\limits_{i}(\a_{i}\vee\diamondsuit \de)\wedge
\be)\Leftrightarrow
\diamondsuit(\underbrace{\bigwedge\limits_{i}\a_{i}\wedge\be}\limits_{\de}),\ \ \ 
\vdash_{G_{Grz}^+}\emptyset|\diamondsuit(\bigwedge\limits_{i}(\a_{i}\vee\diamondsuit \de)\wedge
\be)\Leftrightarrow
\diamondsuit(\underbrace{\bigwedge\limits_{i}\a_{i}\wedge\be}\limits_{\de}).
$$
\end{lem}

\begin{proofof}{Lemma \ref{lem:k4trick}}
The corresponding two implications (same in both cases) are easily seen to hold on transitive models. The claim now follows from the completeness results of Lemma \ref{lem:cutG}. It is however also not hard to write down proofs in $G_{GL}$ and $G_{Grz}^+$, which we omit for space reasons. 
\end{proofof}

\section{Uniform interpolation}

\subsection{Uniform interpolation in GL}
We will prove the uniform interpolation by constructing, for each formula $\a$, the pre-interpolant $\forall p(\a)$. The post-interpolant can be defined by $\exists(\a)=\neg\forall p(\neg\a)$.
The construction is based on a proof-search for the sequent $\emptyset\Rightarrow\a$. To make it work we need to define interpolants for sequents instead of formulas. The uniform interpolation is then obtained via 
$$
\forall p(\a)=\forall p(\emptyset;\a).
$$
 
\begin{thm}\label{thm:GL}
Let $\Gamma,\Delta$ be finite multisets of formulas. For every
propositional variable $p$ there exists a~formula $\forall p(\Gamma;\Delta)$ such that:
\begin{itemize}
\item[(i)]
 $$ Var(\forall p(\Gamma;\Delta))\subseteq Var(\Gamma,\Delta)\backslash \{p\}$$
\item[(ii)]
$$
\vdash_{G_{GL}} \Gamma,\forall p(\Gamma;\Delta)\Rightarrow\Delta
$$
\item[(iii)] moreover let $\Phi ,\Psi$ be multisets of formulas
not containing $p$ and
$$
\vdash_{G_{GL}}\Phi,\Gamma\Rightarrow\Delta,\Psi.
$$ 
Then
$$
 \vdash_{G_{GL}}\Phi\Rightarrow\forall p(\Gamma ; \Delta),\Psi.
$$
\end{itemize}
\end{thm}

\begin{proofof}{Theorem \ref{thm:GL}}
In the following construction of the interpolant, it is instructive to imagine that with the formula $\forall p(\Gamma;\Delta)$ we are describing (a relevant part of) the proof-search tree for the sequent $(\Phi,\Gamma\Rightarrow\Delta,\Psi)$ for any context $\Phi,\Psi$ not containing $p$, namely the part only depending on $\Gamma,\Delta$. This description has to be finite.
The interpolant is defined recursively, closely following the proof-search strategy: for a non-critical sequent we simply use its closure, while for a critical sequent we apply a matching argument, similar to the strategy used by Pitts \cite{Pit92}.
We start with a definition of the formula $\forall p(\Gamma;\Delta)$, then we prove that the definition terminates, and proceed with proving it satisfies items (i)-(iii) of the Theorem \ref{thm:GL}.

\medskip\par\noindent
\textbf{Definition of the interpolant.} We describe the construction of the interpolant recursively. 
The~formula $\forall p(\Gamma ; \Delta)$ is for a noncritical $(\Gamma ; \Delta)$ defined by
\begin{equation}\label{eq:defUI}
\forall p(\Gamma ; \Delta) = \bigwedge\limits_{(\Gamma_i\Rightarrow\Delta_i)\in
Cl(\Gamma;\Delta)}\forall p(\Gamma_i;\Delta_i)
\end{equation}
The~recursive steps for $\Gamma\Rightarrow\Delta$ being a~critical sequent of the form $(\Box\Gamma',\Pi;\Box\Delta',\Lambda)$, with $\Pi,\Lambda$ atomic, are given below in Table \ref{tab:GL}. The~first line of Table \ref{tab:GL} corresponds to some of the~cases when the~critical sequent is provable - it is either an~axiom or the~diagonal formula is already in the~antecedent (here we are using the~loop preventing mechanism from the~termination argument in Lemma \ref{lem:GLterm}).
The line 2 of Table \ref{tab:GL} corresponds to a critical step in a proof-search, the corresponding disjunction covering:
\begin{itemize}
\item[-] propositional variables from multisets $\Pi,\Lambda$, where all $q,r\neq p$,
\item[-] all the~possibilities of a~$\Box_{GL}$ inference with the~principal formula from $\Box\Delta'$,
\item[-] and, by the diamond formula $\diamondsuit\bigwedge N(\Box\Gamma',\Gamma ';\emptyset)$, which we will define below, also the~possibility of
a~$\Box_{GL}$ inference with the~principal formula not from $\Box\Delta'$ (i.e. from a context not containing $p$). Morally, we should include $\diamondsuit\forall p(\Box\Gamma',\Gamma ';\emptyset)$ instead, but then the definition would not terminate. This is the trick we describe below in Remark \ref{rem:trick}.
\end{itemize}

\begin{table}[!htp]
{\renewcommand{\arraystretch}{1.5}
\begin{tabular}{|c|c|c|}
\hline
  & $\Box\Gamma',\Pi;\Box \Delta',\Lambda$ matches & $\forall p(\Box\Gamma',\Pi;\Box\Delta',\Lambda)$ equals\\
\hline
1& if $p\in\Pi\cap\Lambda$ & \\
 & or $\bot\in\Pi$ & \\
  & or $\Gamma'\cap\Delta'\neq\emptyset$ & $\top$\\
\hline
2 & otherwise & $\bigvee\limits_{q\in\Lambda}q\vee\bigvee\limits_{r\in\Pi}\neg r$\\
  & (here all $q,r\neq p$) & $\bigvee\limits_{\be\in\Delta '}\Box\forall p(\Box\Gamma ',\Gamma ',\Box\be;\be)\vee$\\
  & & $\diamondsuit\bigwedge N(\Box\Gamma',\Gamma';\emptyset)$\\
 \hline
\end{tabular}
\caption{ }\label{tab:GL}
}
\end{table}
For a sequent of the form $(\Box\Gamma',\Gamma';\emptyset)$, a set of formulas $N(\Box\Gamma',\Gamma';\emptyset)$ is defined as the smallest set given by Table \ref{tab:NGL}:

\begin{table}[!htp]
{\renewcommand{\arraystretch}{1.5}
\begin{tabular}{|c|c|c|}
\hline
& $(\Box\Sigma,\Upsilon;\Box\Omega,\Theta)\in Cl(\Box\Gamma',\Gamma ';\emptyset)$ matches & $N(\Box\Gamma',\Gamma';\emptyset)$ contains\\
 \hline
1 & $\Sigma^{\circ}\supset\Gamma'^{\circ}$ & \\
 & or $p\in \Upsilon\cap\Theta$ or $\Sigma\cap\Omega\neq\emptyset $& \\
  & or $\bot\in \Upsilon$ & $\forall p(\Box\Sigma,\Upsilon;\Box\Omega,\Theta)$\\
\hline
2 & otherwise & $\bigvee\limits_{q\in\Theta}q\vee\bigvee\limits_{r\in\Upsilon}\neg
     r$ \\
& (here all $q,r\neq p$)  & $\vee\bigvee\limits_{\be\in\Omega}\Box \forall p(\Box\Sigma,\Sigma,\Box\be;\be)$ \\
\hline
\end{tabular}
\caption{ }\label{tab:NGL}
}
\end{table}

In the first line of Table \ref{tab:NGL}, we use the fact that, given $\Sigma^{\circ}\supset\Gamma'^{\circ}$, the sequent $(\Box\Sigma,\Upsilon;\Box\Omega,\Theta)$ is strictly simpler then $(\Box\Gamma',\Gamma ';\emptyset)$ in terms of the measure we will use below to prove termination of the definition, and therefore it is safe to recursively call the procedure. In the remaining cases the sequent is provable and the value in the second column therefore equals $\top$. The second line of Table \ref{tab:NGL} covers the case when $\Sigma^{\circ}=\Gamma'^{\circ}$, and resembles the line 2 of Table \ref{tab:GL}, only to ensure termination, we have omitted the diamond-part of the disjunction.  

\medskip\par\noindent
\textbf{Termination.}
Let us see that the~definition terminates. The~argument is similar to that we have used to prove termination of the~calculus $G_{GL}$ in \ref{lem:GLterm}. Consider a~run of the~procedure for $\forall p(\Phi;\Psi)$ and let $n$
 be the~number of boxed subformulas occurring in $(\Phi;\Psi)$, which bounds the maximal number of critical steps occurring along a~branch in the~tree corresponding to the~run of the~procedure. This is crucial since it enables us to consider an~upper bound of the~weight of an~argument of $\forall p$ occurring during this~run of the~procedure. Put $c=4^{n}w(\Phi;\Psi)$, i.e. an~upper bound of the~weight of an~argument of $\forall p$ occurring during the~run of the~procedure for $(\Phi;\Psi)$. Here, in contrast to the~termination argument for the~calculus $G_{GL}$, we need $4^{n}$ since the~weight of a~recursively called argument of $\forall p$ can increase more. This is because in the~construction of $N(\Box\Gamma',\Gamma';\emptyset)$ we look one level deeper.
Let, for any multiset $\Gamma$, $b(\Gamma)$ be the~number of boxed formulas in $\Gamma$ counted as a~set. For a~$\forall p$ argument $(\Gamma;\Delta)$ occurring during the construction, consider an~ordered pair
$\langle c - b(\Gamma), w(\Gamma ,\Delta) \rangle$. This
measure decreases in every recursive step of the~procedure in terms of the~lexicographical ordering:
\begin{itemize}
\item[-] It is obvious that, for each noncritical sequent
$(\Gamma'\Rightarrow\Delta')\in Cl(\Gamma;\Delta)$,
$w(\Gamma',\Delta')<w(\Gamma,\Delta)$ and that $b$ does not decrease.
\item[-] for a~critical argument $(\Box\Gamma',\Pi;\Box\Delta',\Lambda)$, i.e., line 2 of Table \ref{tab:GL}, and whole Table \ref{tab:NGL} constructing $N(\Box\Gamma',\Gamma';\emptyset)$. For all the~three recursively called arguments $b$ increases, thus $c - b$, and therefore the~whole measure, decreases.
\end{itemize}

\begin{rem}\label{rem:trick}[termination trick]
To retain termination of the definition, we cannot replace 
$$\diamondsuit\bigwedge N(\Box\Gamma',\Gamma';\emptyset)$$ 
in the line 2 of Table \ref{tab:GL} with $\diamondsuit\forall p(\Box\Gamma',\Gamma';\emptyset)$, which in fact seems to be needed to prove the part (iii) of the theorem. The reason is that its recursively called argument need not be simpler then the sequent in question. However, in the prove of (iii) $\diamondsuit\forall p(\Box\Gamma',\Gamma';\emptyset)$ can and will be used, because, as we show next, it is the case that
\begin{equation}\label{eq:trick}
\vdash_{G_{GL}}\diamondsuit\bigwedge N(\Box\Gamma',\Gamma';\emptyset)\Leftrightarrow\diamondsuit\forall p(\Box\Gamma',\Gamma';\emptyset).
\end{equation}
\end{rem}

\par\noindent
To see this, we use the fixed point observation we made earlier in Lemma \ref{lem:k4trick}. Consider sequents $(\Box\Sigma,\Upsilon;\Box\Omega,\Theta)$ in the closure of $(\Box\Gamma',\Gamma';\emptyset)$, and refer by $S$  to sequents with $\Sigma^{\circ}=\Gamma'^{\circ}$, and by $S'$ to sequents in the closure with $\Sigma^{\circ}\supset\Gamma'^{\circ}$, i.e. strictly simpler then $(\Box\Gamma',\Gamma';\emptyset)$.
Since 
$$\forall p(\Box\Gamma',\Gamma';\emptyset)=\bigwedge\limits_{S}\forall p(S)\wedge\bigwedge\limits_{S'}\forall p(S'),
$$
and for each $S = (\Box\Sigma,\Upsilon;\Box\Omega,\Theta)$ with $\Sigma^{\circ}=\Gamma'^{\circ}$ we obtain by the line 2 of Table \ref{tab:GL}
$$
\forall p(S) = \bigvee\limits_{q\in\Theta}q\vee\bigvee\limits_{r\in\Upsilon}\neg r\vee\bigvee\limits_{\be\in\Omega}\Box\forall p(\Box\Sigma,\Sigma,\Box\be;\be)\vee\diamondsuit\bigwedge N(\Box\Sigma,\Sigma;\emptyset),
$$
and using the definition of $N(\Box\Gamma',\Gamma';\emptyset)$ in the Table \ref{tab:NGL}, the above equivalence (\ref{eq:trick}) becomes the following:
\medskip
$$
\vdash_{G_{GL}}\diamondsuit(\bigwedge\limits_{S}
 (\bigvee\limits_{q\in\Theta}q\vee\bigvee\limits_{r\in\Upsilon}\neg r\vee\bigvee\limits_{\be\in\Omega}\Box\forall p(\Box\Sigma,\Sigma,\Box
 \be;\be))\wedge\bigwedge\limits_{S'}\forall p(S')
 $$
 $$
 \Leftrightarrow\ \diamondsuit(\bigwedge\limits_{S}(\bigvee\limits_{q\in\Theta}q\vee\bigvee\limits_{r\in\Upsilon}\neg r\vee\bigvee\limits_{\be\in\Omega}\Box\forall p(\Box\Sigma,\Sigma,\Box\be;\be)\vee\diamondsuit\bigwedge N(\Box\Sigma,\Sigma;\emptyset))\wedge\bigwedge\limits_{S'}\forall p(S')).
$$

\medskip\par\noindent
Observe, that $N(\Box\Sigma,\Sigma;\emptyset)$ is equivalent to $N(\Box\Gamma',\Gamma';\emptyset)$ by $\Sigma^{\circ}=\Gamma'^{\circ}$. Therefore the result now follows by Lemma \ref{lem:k4trick}, instantiated with 
\begin{align*}
\a_S & = & \bigvee\limits_{q\in\Theta}q\vee\bigvee\limits_{r\in\Upsilon}\neg r\vee\bigvee\limits_{\be\in\Omega}\Box\forall p(\Box\Sigma,\Sigma,\Box\be;\be)\\
\be & = & \bigwedge\limits_{S'}\forall p(S')\\
\de & = & \bigwedge N(\Box\Gamma',\Gamma';\emptyset)
\end{align*}
We have thus established the termination of the definition of the uniform interpolants. Now we proceed in proving the three items of Theorem \ref{thm:GL}.
\medskip\par\noindent
\textbf{(i).} The item (i) follows easily by induction on $\Gamma,\Delta$, because we never use $p$ during the~definition of the~formula $\forall p(\Gamma;\Delta)$.

\medskip\par\noindent
\textbf{(ii).} We proceed by induction on the~complexity of $\Gamma,\Delta$ given by the~measure function defined above, and prove that $\vdash_{G_{GL}}\Gamma,\forall p(\Gamma;\Delta)\Rightarrow\Delta$. 

First, let $(\Gamma\Rightarrow\Delta)$ be a~noncritical sequent. Then sequents $(\Gamma_{i}\Rightarrow\Delta_{i})\in Cl(\Gamma;\Delta)$ are of lower complexity and by the~induction hypotheses
$\vdash_{G_{GL}}\Gamma_{i},\forall p(\Gamma_{i};\Delta_{i})\Rightarrow\Delta_{i}$
for each i. Then by admissibility of weakening and by Lemma \ref{lem:cl}
$$
\vdash_{G_{GL}}\Gamma,\forall p(\Gamma_{1};\Delta_{1}),\ldots,\forall p(\Gamma_{k};\Delta_{k})\Rightarrow\Delta,
$$
therefore by a $\wedge$-l inference
$$
\vdash_{G_{GL}}\Gamma,\bigwedge\limits_{(\Gamma_{i}\Rightarrow\Delta_{i})\in
Cl(\Gamma;\Delta)} \forall p(\Gamma_{i};\Delta_{i})\Rightarrow\Delta,
$$ 
which is by (\ref{eq:defUI})
$$
\vdash_{G_{GL}}\Gamma, \forall p(\Gamma;\Delta)\Rightarrow\Delta.
$$
Second, let $(\Gamma\Rightarrow\Delta)$ be a~critical sequent. 
If $(\Gamma\Rightarrow\Delta)$ is a~critical sequent matching
the~line 1 of Table \ref{tab:GL}, then (ii) is an~axiom or a provable sequent.
Let $(\Gamma\Rightarrow\Delta)$ be a~critical sequent matching
the~line 2 of Table \ref{tab:GL}. We prove 
$$
\vdash_{G_{GL}}\Pi,\de,\Box\Gamma'\Rightarrow\Box\Delta',\Lambda
$$ 
for each disjunct $\de$ used in the line 2 of Table \ref{tab:GL} to define the interpolant:
\begin{itemize}
\item[-] For each $r\in\Pi$ obviously $\vdash_{G_{GL}}\Pi,\neg
r,\Box\Gamma'\Rightarrow\Box\Delta',\Lambda$, therefore \\ $\vdash_{G_{GL}}\Pi,\bigvee\limits_{r\in\Pi}\neg
r,\Box\Gamma'\Rightarrow\Box\Delta',\Lambda$.

\item[-] for each $q\in\Lambda$ obviously $\vdash_{G_{GL}}\Pi
,q,\Box\Gamma'\Rightarrow\Box\Delta',\Lambda$, therefore \\ $\vdash_{G_{GL}}\Pi,\bigvee\limits_{q\in\Lambda} q,\Box\Gamma'\Rightarrow\Box\Delta',\Lambda$ 

\item[-] For each $\be\in\Delta '$, $\vdash_{G_{GL}}\Box\Gamma
',\Gamma ', \Box\be,\forall p(\Box\Gamma',\Gamma ',\Box\be;\be)
\Rightarrow\be$ by the~induction hypothesis, which gives
$\vdash_{G_{GL}}\Box\Gamma',\Pi,\Box\forall p(\Box\Gamma',\Gamma
',\Box\be;\be)\Rightarrow\Box\Delta',\Lambda$ by
a~$\Box_{L}$ inference.
\end{itemize}
\medskip
It remains to be proved that 
$$
\vdash_{G_{GL}}\Box\Gamma',\Pi,\diamondsuit\bigwedge N(\Box\Gamma',\Gamma';\emptyset)\Rightarrow\Box\Delta',\Lambda.
$$
For each $(\Box\Sigma,\Upsilon\Rightarrow\Box\Omega,\Theta)\in
Cl(\Box\Gamma',\Gamma';\emptyset)$ from the first line of Table \ref{tab:NGL}, we know by the~induction hypotheses that
$$
\vdash_{G_{GL}}\Box\Sigma,\Upsilon,\forall p(\Box\Sigma,\Upsilon;\Box\Omega,\Theta)\Rightarrow\Box\Omega,\Theta
.$$
For each $(\Box\Sigma,\Upsilon\Rightarrow\Box\Omega,\Theta)\in
Cl(\Box\Gamma ',\Gamma';\emptyset)$ from the second line of Table \ref{tab:NGL}
we have the following:
\begin{itemize}
\item[-]
$
\vdash_{G_{GL}}\Box\Sigma,\Upsilon,\bigvee\limits_{q\in\Theta}q\Rightarrow\Box\Omega,\Theta.
$
\item[-] 
$
\vdash_{G_{GL}}\Box\Sigma,\Upsilon,\bigvee\limits_{\neg
r\in\Upsilon}\neg r\Rightarrow\Box\Omega,\Theta.
$
\item[-] for each
$\be\in\Omega$ by the induction hypotheses
$$
\vdash_{G_{GL}}\Box\Sigma,\Sigma,\Box\be,\forall p(\Box\Sigma,\Sigma,\Box\be;\be)\Rightarrow\be$$ and by weakening
and a~$\Box_{GL}$ inference 
$$
\vdash_{G_{GL}}\Box\Sigma,\Box\forall p(\Box\Sigma,\Sigma,\Box\be;\be)\Rightarrow\Box\be.
$$
\end{itemize}
Together this yields, using $\vee$-l inferences,
$$
\vdash_{G_{GL}}\Box\Sigma,\Upsilon, \bigvee\limits_{q\in\Theta}\vee\bigvee\limits_{\neg r\in\Upsilon}
\vee\bigvee\limits_{\be\in\Omega}\Box\forall p(\Box\Sigma,\Sigma,\Box\be;\be)\Rightarrow\Box\Omega,\Theta.
$$

\medskip\par\noindent
Therefore, for each $(\Box\Sigma,\Upsilon\Rightarrow\Box\Omega,\Theta)\in
Cl(\Box\Gamma ',\Gamma ';\emptyset)$, we obtain, using weakening and $\wedge$-r inferences,
$$
\vdash_{G_{GL}}\Box\Sigma,\Upsilon,\bigwedge N(\Box\Gamma',\Gamma';\emptyset)\Rightarrow\Box\Omega,\Theta.
$$
Now, by Lemma \ref{lem:cl}, 
$$
\vdash_{G_{GL}}\Box\Gamma ',\Gamma',\bigwedge N(\Box\Gamma',\Gamma';\emptyset)\Rightarrow\emptyset
.$$ 
By negation and weakening inferences 
$$
\vdash_{G_{GL}}\Box\Gamma',\Gamma',\Box\neg\bigwedge N(\Box\Gamma',\Gamma';\emptyset)\Rightarrow\neg\bigwedge N(\Box\Gamma',\Gamma';\emptyset)
$$ 
and by a~$\Box_{GL}$ inference
$$
\vdash_{G_{GL}}\Box\Gamma',\Pi\Rightarrow\Box\neg
\bigwedge N(\Box\Gamma',\Gamma';\emptyset),\Box\Delta',\Lambda.
$$ Now, using a~negation inference again, we
obtain 
$$
\vdash_{G_{GL}}\Box\Gamma',\Pi,\diamondsuit\bigwedge N(\Box\Gamma',\Gamma';\emptyset)\Rightarrow\Box\Delta',\Lambda.
$$

\medskip\par\noindent
Putting finally all the above disjuncts together then yields, using $\vee$-l inferences,
$$
\vdash_{G_{GL}}\Pi,\Box\Gamma',\bigvee\limits_{q\in\Lambda}q
\bigvee\limits_{r\in\Pi}\neg r\bigvee\limits_{\be\in\Delta'}\Box
\forall p(\Box\Gamma ',\Gamma ',\Box\be;\be)\vee\diamondsuit\bigwedge N(\Box\Gamma',\Gamma';\emptyset)\Rightarrow\Box\Delta',\Lambda,
$$ 
that is, by the~line 2 of Table \ref{tab:GL},
$$
\vdash_{G_{GL}}\Pi,\Box\Gamma',\forall p(\Pi,\Box\Gamma';\Box\Delta',\Lambda)\Rightarrow\Box\Delta',\Lambda.
$$

\medskip\par\noindent
\textbf{(iii).} We proceed by induction on the~height of a~proof of
$(\Phi,\Gamma\Rightarrow\Delta,\Psi)$, and by sub-induction on the measure of the sequent $(\Gamma;\Delta)$ used to show termination of the definition. We show that $\vdash_{G_{GL}}\Phi\Rightarrow\forall p(\Gamma;\Delta),\Psi$.

\medskip\par\noindent
First, consider $(\Phi,\Gamma\Rightarrow\Delta,\Psi)$ is an axiom. The following cases apply:
\begin{itemize}
\item[-] $\bot$ is principal and $\bot\in\Phi$, then (iii) is an axiom.
\item[-] $\bot$ is principal and $\bot\in\Gamma$, then $\forall p(\Gamma;\Delta)=\top$ and $\Phi\Rightarrow\top,\Psi$ is provable. 
\item[-] $p$ is principal, i.e. $p\in\Gamma\cap\Delta$ and $\forall p(\Gamma;\Delta)=\top$ and $\Phi\Rightarrow\top,\Psi$ is provable. 
\item[-] $q\neq p$ is principal, and $q\in\Phi\cap\Psi$. Then $\Phi\Rightarrow\forall p(\Gamma;\Delta),\Psi$ is an axiom.
\item[-] $q\neq p$ is principal, and $q\in\Phi\cap\Delta$. Then $\vdash_{G_{GL}}q\Rightarrow\forall p(\Gamma;\Delta)$ by the line 1 of Table \ref{tab:GL}, and we obtain the result by weakening.
\item[-] $q\neq p$ is principal, and $q\in\Gamma\cap\Psi$. Then $\vdash_{G_{GL}}\neg q\Rightarrow\forall p(\Gamma;\Delta)$ by the line 1 of Table \ref{tab:GL}, and $\vdash_{G_{GL}}\emptyset\Rightarrow\forall p(\Gamma;\Delta), q$ by $\neg$-l invertibility, and we obtain the result by weakening.
\item[-] $q\neq p$ is principal, and $q\in\Gamma\cap\Delta$. Then $\vdash_{G_{GL}}q\vee\neg q\Rightarrow\forall p(\Gamma;\Delta)$ by the line 1 of the table, and therefore $\vdash_{G_{GL}}\emptyset\Rightarrow\forall p(\Gamma;\Delta)$, and we obtain the result by weakening.
\end{itemize}

\medskip\par\noindent
Consider then $(\Phi,\Gamma\Rightarrow\Delta,\Psi)$ is not an axiom. We distinguish two main cases: Consider first $(\Gamma;\Delta)$ is a noncritical sequent. Then all $(\Gamma'\Rightarrow\Delta')\in Cl(\Gamma;\Delta)$ are strictly simpler in terms of the measure, and for all of them 
we have ${\vdash_{G_{GL}}\Phi,\Gamma'\Rightarrow\Delta',\Psi}$ by Lemma
\ref{lem:cl}. Then, using the induction hypothesis and (\ref{eq:defUI}),  the~following are equivalent:
$$
\vdash_{G_{GL}}\Phi\Rightarrow\forall p(\Gamma';\Delta'),\Psi\quad \mbox{for
all}\quad (\Gamma'\Rightarrow\Delta')\in Cl(\Gamma;\Delta)
$$
$$
\vdash_{G_{GL}}\Phi\Rightarrow
\bigwedge\limits_{(\Gamma'\Rightarrow\Delta')\in Cl(\Gamma;\Delta)}
\forall p(\Gamma';\Delta'),\Psi
$$
$$
\vdash_{G_{GL}}\Phi\Rightarrow\forall p(\Gamma;\Delta),\Psi.
$$

\medskip\par\noindent
Consider $(\Gamma;\Delta)$ is a critical sequent and the last inference is an instance of an invertible rule. Then the principal formula of the inference is in $\Phi,\Psi$. We apply the induction hypothesis to the premise of the last inference, and then the invertible rule in question again.

\medskip\par\noindent
Finally assume that $(\Gamma;\Delta)$ is a critical sequent and the~last inference is a~$\Box_{L}$ inference:

\medskip\par\noindent
Consider first that the~principal formula $\Box\a\in\Psi$, in particular, $\a$ doesn't contain $p$. Then the~proof ends with the step:
\begin{prooftree}
\AxiomC{$\Box\Phi',\Box\Gamma ',\Phi',\Gamma',\Box\a\Rightarrow\a$}
\RL{$\Box_{L}$}
\UnaryInfC{$\Box\Pi',\Box\Gamma',\Pi'',\Gamma''\Rightarrow\Box\a,\Psi',\Delta$}
\end{prooftree}
where $\Box\Phi',\Phi''$ is $\Phi$; $\Box\Gamma',\Gamma''$ is
$\Gamma$; and $\Box\a,\Psi'$ is $\Psi$. Consider
$\Box\Gamma'\cap\Delta =\emptyset$ (otherwise
the line 1 of Table \ref{tab:GL} applies and $\forall p(\Gamma;\Delta)=\top$, and therefore (iii) holds). So we can use
the~line 2 of Table \ref{tab:GL}. Then the~induction hypothesis gives
$$
\vdash_{G_{GL}}\Phi',\Box\a\Rightarrow\forall p(\Box\Gamma',\Gamma';\emptyset),\a
$$ 
and by a~$\neg$-l inference we obtain 
$$
\vdash_{G_{GL}}\Phi',\Box\a,\neg
\forall p(\Box\Gamma',\Gamma';\emptyset)\Rightarrow\a.
$$ 
Now, by
a~$\Box_L$ and a~negation inferences, we obtain
$$
\vdash_{G_{GL}}\Box\Phi',\Phi''\Rightarrow\diamondsuit\forall p(\Box\Gamma',\Gamma';\emptyset),
\Box\a,\Psi'.
$$ 
By the~line 2 of Table \ref{tab:GL}, invertibility of
the~$\vee$-l rule, and by (\ref{eq:trick}) we have
$$
\vdash_{G_{GL}}\diamondsuit
\forall p(\Box\Gamma',\Gamma';\emptyset)\Rightarrow
\forall p(\Box\Gamma',\Gamma'';\Delta).
$$ 
The two sequents above yield (iii) by cut admissibility.

\medskip\par\noindent
Consider the~principal formula $\Box\a\in\Delta$. Again, consider
$\Box\Gamma'\cap\Delta =\emptyset$ so we can use the~line 2 of Table \ref{tab:GL}.
Then the~proof ends with:
\begin{prooftree}
\AxiomC{$\Box\Phi',\Box\Gamma',\Phi',\Gamma',\Box\a\Rightarrow\a$}
\RL{$\Box_L$}
\UnaryInfC{$\Box\Phi',\Box\Gamma',\Phi'',\Gamma''\Rightarrow\Box\a
,\Delta ',\Psi$}
\end{prooftree}
where $\Box\Phi',\Phi''$ is $\Phi$; $\Box\Gamma',\Gamma''$ is
$\Gamma$; and $\Box\a,\Delta'$ is $\Delta$. Now the~induction hypothesis gives 
$$
\vdash_{G_{GL}}\Box\Phi',\Phi'
\Rightarrow\forall p(\Box\Gamma',\Gamma ',\Box\a;\a),
$$ 
and by weakening and a~$\Box_L$ inference we obtain
$$
\vdash_{G_{GL}}\Box\Phi',\Phi''\Rightarrow\Box\forall p(\Box\Gamma',\Gamma ',\Box\a;\a),\Psi.
$$
The~line 2 of Table \ref{tab:GL} and invertibility of the~$\vee$-l rule yields
$$
\vdash_{G_{GL}}\Box\forall p(\Box\Gamma',\Gamma',\Box\a;\a)\Rightarrow \forall p(\Box\Gamma',\Gamma'';\Box\a,\Delta').
$$ 
Finally, we obtain (iii) by cut admissibility.
\end{proofof}

\begin{rem}[\textbf{Constructivity of the proof}] We have given a construction of uniform interpolant which is effective and implementable. However, since we argued semantically to claim cut-free completeness of the calculus, reader might object that our proof of the uniform interpolation theorem is not fully constructive. To this point we say the following: one can look at the cut-elimination proof in \cite{GR08,GR12} and prove constructively that the two calculi are equivalent. Another way is to use the proof-search procedure described in subsection \ref{sub:termGL} and prove completeness via decidability. The point is that an unsuccessful proof-search tree can be used to construct a counterexample to a given sequent, in spirit of the proof contained in \cite{Sve00}. We have not included such an argument here mainly for space reasons and because it is not essential to understand the proof of uniform interpolation.
\end{rem}
\subsubsection{Fixed points}\label{subsec:fixedp}

Uniform interpolation theorem for \textbf{GL} entails Sambin's
and de Jongh's fixed point theorem. Our proof then presents
an~alternative constructive proof of the~fixed point theorem:

\begin{thm}\label{thm:fp} Fixed point theorem: Suppose $p$ is modalized in $\be$
(i.e., any occurrence of $p$ is in the~scope of a~$\Box$). Then we
can find a~formula $\g$ in the~variables of $\be$ without $p$ such that
$$\vdash_{GL}\g\leftrightarrow \be(\g).$$
\end{thm}

Already Craig interpolation entails fixed point theorem: a~fixed point of a~formula $\be$ is
an~interpolant of a~sequent expressing
the~uniqueness of the~fixed point
$$
\boxdot(p\leftrightarrow \be(p))\wedge\boxdot (q\leftrightarrow \be(q))
\Rightarrow p\leftrightarrow q ,
$$ 
which is provable in \textbf{GL}
- proofs of this fact in \cite{SV82} and \cite{Boo93} are easily adaptable to our variant of the calculus. However,
to construct the fixed point using this method requires to have an~actual proof of the~sequent expressing the~uniqueness.

Direct proofs of fixed point theorem were given by Sambin
\cite{Sam76}, Sambin and Valentini in \cite{SV82} (a construction of
explicit fixed points which is effective and implementable),
Smory\`{n}ski \cite{Smo78} from Beth's definability property,
Reidhar-Olson \cite{LRO90}, Gleit and Goldfarb \cite{GG90}. A~proof
from Beth's property can be found also in Kracht's book
\cite{Kra99}, for three different proofs see Boolos' book
\cite{Boo93}. A~different and effective constructive proof of fixed point theorem is the~one by Sambin and Valentini in \cite{SV82}. 
We present a proof of fixed point theorem based on uniform interpolation, it is an~effective proof alternative to those above. We learnt this simple argument from Albert Visser, and we found it an interesting application of the uniform interpolation theorem.

\begin{proofof}{Theorem \ref{thm:fp}}
Let us consider a~formula $\be(p,\bar{q})$ with $p$ modalized in $\be$.
The~fixed point of $\be$ then would be the~simulation of 
$$
\exists
p(\Box (p\leftrightarrow \be(p))\wedge \be(p))
$$ 
or, equivalently, of
$$
\forall r(\Box (r\leftrightarrow \be(r))\rightarrow \be(r)).
$$ 
Let us denote
them $\g_{1}$ and $\g_{2}$ and observe they are both interpolants of the sequent
$$
(\Box (p\leftrightarrow \be(p))\wedge \be(p)\Rightarrow \Box
(r\leftrightarrow \be(r))\rightarrow \be(r))
$$ 
and that neither of them contains $p,r$. We show that any of them is the~fixed point of $\be(p)$
and that they are indeed equivalent. To keep readability we just sketch the proofs in $G_{GL}$ below. 

First we show that $(\Box (p\leftrightarrow \be(p))\wedge
\be(p)\Rightarrow \Box (r\leftrightarrow \be(r))\rightarrow \be(r))$ is
provable from the~uniqueness statement:
\begin{prooftree}
\AxiomC{$\boxdot(p\leftrightarrow \be(p))\wedge\boxdot
(r\leftrightarrow \be(r))\Rightarrow p\leftrightarrow
r$}\AxiomC{$\vdots$}\UnaryInfC{$p\leftrightarrow r,p\leftrightarrow
\be(p),r\leftrightarrow \be(r),\be(p)\Rightarrow
\be(r)$}\RightLabel{cut}\BinaryInfC{$\Box(p\leftrightarrow
\be(p)),\Box(r\leftrightarrow \be(r)),\be(p)\Rightarrow
\be(r)$}\UnaryInfC{$(\Box (p\leftrightarrow \be(p)), \be(p)\Rightarrow
\Box (r\leftrightarrow \be(r))\rightarrow \be(r))$}\UnaryInfC{$(\Box
(p\leftrightarrow \be(p))\wedge \be(p)\Rightarrow \Box (r\leftrightarrow
\be(r))\rightarrow \be(r))$}
\end{prooftree}

Now let us see that any of $\g_{i}$ is a~fixed point and thus, by
the~uniqueness, $\g_{1}\leftrightarrow \g_{2}$. First observe, that whenever $(\Gamma(p)\Rightarrow\Delta(p))$ is
provable, $(\Gamma[p/\a]\Rightarrow\Delta[p/\a])$ where we substitute
$\a$ for $p$ is provable as well (we substitute everywhere in
the~proof, to treat $\Box_{GL}$ inferences can require some
admissible weakenings, and we add proofs of sequents $(\Gamma,
\a\Rightarrow\Delta,\a)$ in place of axioms with $p$ principal). The~label
"subst." in the~following proof-tree refers to such a~substitution,
the~label "inv." refers to invertibility of a~rule:
\begin{prooftree}
\AxiomC{$\Box (p\leftrightarrow \be(p))\wedge \be(p)\Rightarrow\g_{i}$}
\RightLabel{subst.}
\UnaryInfC{$\Box (\g_{i}\leftrightarrow\be(\g_{i}))\wedge\be(\g_{i})\Rightarrow\g_{i}$}
\RightLabel{inv.}
\UnaryInfC{$\Box (\g_{i}\leftrightarrow\be(\g_{i})),\be(\g_{i})\Rightarrow\g_{i}$}
\UnaryInfC{$\Box (\g_{i}\leftrightarrow\be(\g_{i}))\Rightarrow \neg\be(\g_{i}),\g_{i}$}
\UnaryInfC{$\Box(\g_{i}\leftrightarrow\be(\g_{i}))\Rightarrow\be(\g_{i})\rightarrow\g_{i}$} 
\AxiomC{$\g_{i}\Rightarrow\Box (r\leftrightarrow\be(r))\rightarrow\be(r)$}
\RightLabel{subst.}
\UnaryInfC{$\g_{i}\Rightarrow\Box(\g_{i}\leftrightarrow\be(\g_{i}))\rightarrow\be(\g_{i})$}
\RightLabel{inv.} 
\UnaryInfC{$\g_{i}\Rightarrow\neg\Box(\g_{i}\leftrightarrow\be(\g_{i})),\be(\g_{i})$} 
\UnaryInfC{$\g_{i},\Box(\g_{i}\leftrightarrow\be(\g_{i}))\Rightarrow\be(\g_{i})$}
\UnaryInfC{$\Box (\g_{i}\leftrightarrow\be(\g_{i}))\Rightarrow \neg\g_{i},\be(\g_{i})$}
\UnaryInfC{$\Box(\g_{i}\leftrightarrow\be(\g_{i}))\Rightarrow\g_{i}\rightarrow\be(\g_{i})$}
\BinaryInfC{$\Box (\g_{i}\leftrightarrow\be(\g_{i}))\Rightarrow\g_{i}\leftrightarrow\be(\g_{i})$}
\RightLabel{$\Box_{GL}$}
\UnaryInfC{$\emptyset\Rightarrow\Box(\g_{i}\leftrightarrow\be(\g_{i}))$}
\end{prooftree}
Now by a~cut
\begin{prooftree}
\AxiomC{$\emptyset\Rightarrow\Box(\g_{i}\leftrightarrow\be(\g_{i}))$}
\AxiomC{$\Box(\g_{i}\leftrightarrow\be(\g_{i}))\Rightarrow\g_{i}\leftrightarrow\be(\g_{i})$}
\RightLabel{cut}
\BinaryInfC{$\emptyset\Rightarrow\g_{i}\leftrightarrow\be(\g_{i})$}
\end{prooftree}

From this proof one can see that already ordinary interpolation does
the~job. The~point of using uniform interpolation here is that we do
not need to have an actual~proof of $(\Box (p\leftrightarrow \be(p))\wedge\be(p)\Rightarrow\Box(r\leftrightarrow\be(r))\rightarrow\be(r))$ to construct a fixed point - we just need to know that the sequent is provable to show that we have indeed constructed a fixed point.

\end{proofof}

\subsection{Uniform interpolation in Grz}
The proof of uniform interpolation in \textbf{Grz} follows the same ideas and is very similar to the previous one, only syntactically a bit more more complicated.

\begin{thm}\label{thm:Grz}
Let $\Gamma,\Delta,\Sigma$ be finite multisets of formulas. For every propositional
variable $p$ there exists a~formula $\forall p(\Box\Sigma |\Gamma;\Delta)$
such that:
\begin{itemize}
\item[(i)]
$$ 
 Var(\forall p(\Box\Sigma |\Gamma ; \Delta))\subseteq Var(\Sigma ,\Gamma,\Delta)\backslash \{p\}
$$
\item[(ii)]
$$
\vdash_{G_{Grz}^{+}}\Box\Sigma|\Gamma ,\forall p(\Box\Sigma |\Gamma ; \Delta)\Rightarrow\Delta 
$$
\item[(iii)] moreover let $\Phi ,\Psi,\Theta$ be multisets of
formulas not containing $p$ and
$$
\vdash_{G_{Grz}^{+}}\Box\Theta,\Box\Sigma|\Phi,\Gamma\Rightarrow\Psi,\Delta.
$$
Then
$$ 
\vdash_{G_{Grz}^{+}}\emptyset|\Box\Theta,\Phi \Rightarrow \forall p(\Box\Sigma|\Gamma;\Delta),\Psi.
$$
\end{itemize}
\end{thm}

\begin{proofof}{Theorem \ref{thm:Grz}}
We start with a definition of the formula $\forall p(\Box\Sigma |\Gamma;\Delta)$, then we prove that the definition terminates, and proceed with proving it satisfies items (i)-(iii) of the Theorem. We remark that the item (iii) is formulated with the third multiset empty because we only have a particular form of cut admissible, see Remark \ref{rem:cut}.

\medskip\par\noindent
\textbf{Definition of the interpolant.} We describe the construction of the interpolant recursively. 
The~formula $\forall p(\Box\Sigma|\Gamma;\Delta)$ is defined by
\begin{equation}\label{eq:defUI2}
\forall p(\Box\Sigma|\Gamma;\Delta) = \bigwedge\limits_{(\Box\Sigma_i|\Gamma_I\Rightarrow\Delta_i)\in
Cl(\Box\Sigma |\Gamma;\Delta)}\forall p(\Box\Sigma_i|\Gamma_i;\Delta_i)
\end{equation}

The~recursive steps for $(\Box\Sigma |\Gamma\Rightarrow\Delta)$ being a~critical sequent of the form $(\Box\Gamma'|\Pi;\Box\Delta',\Lambda)$, with $\Pi,\Lambda$ atomic, are given by the~following table:

\begin{table}[!hp]
{\renewcommand{\arraystretch}{1.5}
\begin{tabular}{|c|c|c|}
\hline
  & $\Box\Gamma'|\Pi;\Box \Delta',\Lambda$ matches & $\forall p(\Box\Gamma'|\Pi;\Box\Delta',\Lambda)$ equals\\
\hline
1& if $p\in\Pi\cap\Lambda$ & \\
& or $\bot\in\Pi$ & \\
 & or $\Gamma'\cap\Delta'\neq\emptyset$ & $\top$\\
\hline
2 & otherwise & $\bigvee\limits_{q\in\Lambda}q\vee\bigvee\limits_{r\in\Pi}\neg r$\\
 & (here all $q,r\neq p$) & $\bigvee\limits_{\be\in\Delta ',D(\be)\notin\Box\Gamma'}\Box\forall p(\Box\Gamma ',D(\be) |\Gamma ';\be)$\\
 & & $\bigvee\limits_{\be\in\Delta ',D(\be)\in\Box\Gamma'}\Box\forall p(\Box\Gamma'|\emptyset;\be))$\\
 &  & $\vee\diamondsuit \bigwedge N(\Box\Gamma'|\Gamma';\emptyset)$\\
 \hline
\end{tabular}
\caption{ }\label{tab:Grz}
}
\end{table}
As in Table \ref{tab:GL} before, the~first line corresponds to some of the~cases when the~critical sequent is provable, and the line 2 corresponds to a critical step, the corresponding disjunction covering
\begin{itemize}
\item[-] propositional variables from multisets $\Pi,\Lambda$, 
\item[-] all the~possibilities of $\Box_{Grz1}^+$ and $\Box_{Grz2}^+$ inferences with the~principal formula from $\Box\Delta'$,
\item[-] and, by the diamond formula $\diamondsuit \bigwedge N(\Box\Gamma'|\Gamma';\emptyset)$ defined below in Table \ref{tab:NGrz}, also the~possibility of
a~$\Box_{Grz1}^+$ or a $\Box_{Grz2}^+$ inference with the~principal formula not from $\Box\Delta'$ (i.e. from a context not containing $p$). For a sequent of the form $(\Box\Gamma'|\Gamma';\emptyset)$, a set of formulas $N(\Box\Gamma'|\Gamma';\emptyset)$ is defined as the smallest set given by the Table \ref{tab:NGrz}.
\end{itemize}

\begin{table}[!htp]
{\renewcommand{\arraystretch}{1.5}
\begin{tabular}{|c|c|c|}
\hline
 & $(\Box\Sigma|\Upsilon;\Box\Omega,\Theta)\in Cl(\Box\Gamma'|\Gamma ';\emptyset)$ matches & $N(\Box\Gamma',\Gamma';\emptyset)$ contains\\
 \hline
 1 & $\Sigma^{\circ}\supset\Gamma'^{\circ}$ & \\
 & or $p\in\Upsilon\cap\Theta$, or $\Sigma\cap\Omega\neq\emptyset$ & \\
 & or $\bot\in\Upsilon$ & $\forall p(\Box\Sigma|\Upsilon;\Box\Omega,\Theta)$\\
\hline
2 & otherwise & $\bigvee\limits_{q\in\Theta}q\vee\bigvee\limits_{r\in\Upsilon}\neg r\vee$\\
 & (here all $q,r\neq p$) & $\bigvee\limits_{\be\in\Omega,D(\be)\notin\Box\Sigma}\Box\forall p(\Box\Sigma,D(\be)|\Sigma;\be)\vee$\\
& & $\bigvee\limits_{\be\in\Omega,D(\be)\in\Box\Sigma}\Box\forall p(\Box\Sigma|\emptyset;\be))$ \\
\hline
\end{tabular}
\caption{ }\label{tab:NGrz}
}
\end{table}

\medskip\par\noindent
\textbf{Termination.} We adopt the~same simplification as we have used proving
termination of the~calculus $G_{Grz^{+}}$ --- we treat the~third multiset as a set (i.e.,
we remove duplicate formulas stored in the set). Consider a~run of
the~procedure for $\forall p(\emptyset|\Phi;\Psi)$.  Let $n$ be
the number of boxed subformulas occurring in $\Phi;\Psi$, which is, as in the~case of
\textbf{GL}, maximal number of critical steps along one branch of
the~corresponding tree.
With each $\forall p$ argument $(\Box\Sigma|\Gamma;\Delta)$ occurring during
the~run of the~procedure, we associate an~ordered pair $\langle n^2-
|\Sigma^{\circ}|,w(\Gamma,\Delta) \rangle$, 
where $n^2$ is an~upper bound of the~number of
formulas stored in $\Box\Sigma$ if we do not duplicate them.
The~measure strictly decreases in each step of the~run of the~procedure in
terms of the~lexicographical ordering:
\begin{itemize}
\item[-]
For a~noncritical argument $(\Box\Sigma|\Gamma;\Delta)$ and for each
$(\Box\Sigma'|\Gamma';\Delta')\in Cl(\Box\Sigma|\Gamma;\Delta)$,
$w(\Gamma',\Delta')<w(\Gamma ,\Delta)$.

\item[-]
For a~critical argument $(\Box\Gamma'|\Pi;\Box \Delta',\Lambda)$
let us see that, in Table \ref{tab:Grz} and Table \ref{tab:NGrz}, for each of the~five recursively called arguments the~measure decreases.
\begin{itemize}
    \item[-] the line 2 in Table \ref{tab:Grz}, $(\Box\Gamma',D(\be)|\Gamma ';\be)$
    where $\be\in\Delta '$ and $D(\be)\notin\Box\Gamma' $: 
    here obviously $|(\Box\Gamma'\cup D(\be))^{\circ}|>|\Box\Gamma
    '^{\circ}|$.

    \item[-] the line 2 in Table \ref{tab:Grz}, $(\Box\Gamma'|\emptyset;\be)$ where $\be\in\Delta '$ and $D(\be)\in\Box\Gamma'$:\\ 
    in this case, $w(\emptyset,\be)<w(\Pi,\Box\Delta',\Lambda)$.
    
    \item[-] the first line in Table \ref{tab:NGrz},  $(\Box\Sigma|\Upsilon;\Box\Omega,\Theta)$ where ${(\Box\Sigma|\Upsilon\Rightarrow\Box\Omega,\Theta)\in Cl(\Box\Gamma '|\Gamma
    ' \emptyset)}$ and $\Sigma^{\circ}\supset\Gamma'^{\circ}$: 
    Since $\Sigma^{\circ}\supset\Gamma'^{\circ}$,
    $|\Sigma^{\circ}|>|\Gamma'^{\circ}|$.

    \item[-] the second line in Table \ref{tab:NGrz}, $(\Box\Sigma,D(\be)|\Sigma;\be)$ where $\be\in\Omega$, 
    $D(\be)\notin\Box\Sigma$, and $(\Box\Sigma|\Upsilon\Rightarrow\Box\Omega,\Theta)\in Cl(\Box\Gamma'|\Gamma ';\emptyset)$ with
    $\Sigma^{\circ}=\Gamma'^{\circ}$:\\
    Since
    $\Sigma^{\circ}=\Gamma'^{\circ}$, also
    $|\Box\Sigma^{\circ}|=|\Box\Gamma'^{\circ}|$.
    Hence $|(\Box\Sigma\cup D(\be))^{\circ}|>|\Box\Gamma'^{\circ}|$.

    \item[-] the second line in Table \ref{tab:NGrz}, $(\Box\Sigma|\emptyset;\be)$ where $\be\in\Omega$, $D(\be)\in\Sigma$,\\  ${(\Box\Sigma|\Upsilon\Rightarrow\Box\Omega,\Theta)\in Cl(\Box\Gamma'|\Gamma';\emptyset)}$ and
    $\Sigma^{\circ}=\Gamma'^{\circ}$:\\
    here $w(\emptyset,\be)<w(\Pi,\Box\Delta',\Lambda)$. 
\end{itemize}
\end{itemize}
\begin{rem}[termination trick] 
Analogously to (\ref{eq:trick}), we want to use $\forall p(\Box\Gamma'|\Gamma';\emptyset)$ in place of $\diamondsuit\bigwedge N(\Box\Gamma'|\Gamma';\emptyset)$ while proving the item (iii) of the theorem. We show next that it is indeed the case that
\begin{equation}\label{eq:trick2}
\vdash_{G_{Grz}^+}\diamondsuit\bigwedge N(\Box\Gamma'|\Gamma';\emptyset)\Leftrightarrow\diamondsuit\forall p(\Box\Gamma'|\Gamma';\emptyset).
\end{equation}
Consider sequents $(\Box\Sigma|\Upsilon;\Box\Omega,\Theta)$ in the closure of $(\Box\Gamma'|\Gamma';\emptyset)$, and refer by $S$  to sequents with $\Sigma^{\circ}=\Gamma'^{\circ}$, and by $S'$ to sequents in the closure with $\Sigma^{\circ}\supset\Gamma'^{\circ}$, i.e. strictly simpler then $(\Box\Gamma'|\Gamma';\emptyset)$.
Since 
$$\forall p(\Box\Gamma'|\Gamma';\emptyset)\equiv\bigwedge\limits_{S}\forall p(S)\wedge\bigwedge\limits_{S'}\forall p(S'),
$$
and for each $S = (\Box\Sigma|\Upsilon;\Box\Omega,\Theta)$ with $\Sigma^{\circ}=\Gamma'^{\circ}$ we obtain $\forall p(S)$ by the line 2 of Table \ref{tab:Grz} to be the following formula:
$$
\bigvee\limits_{q\in\Theta}q\vee\bigvee\limits_{r\in\Upsilon}\neg r\vee\bigvee\limits_{\substack{\be\in\Omega \\D(\be)\notin\Box\Sigma}}\Box\forall p(\Box\Sigma,D(\be)|\Sigma;\be)\vee\bigvee\limits_{\substack{\be\in\Omega \\D(\be)\in\Box\Sigma}}\Box\forall p(\Box\Sigma|\emptyset;\be))\vee\diamondsuit\bigwedge N(\Box\Sigma,\Sigma;\emptyset),
$$
which we can shorten as
$$
\a_S\vee\diamondsuit\bigwedge N(\Box\Sigma,\Sigma;\emptyset).
$$
Now using the definition of $N(\Box\Gamma'|\Gamma';\emptyset)$ in Table \ref{tab:NGrz}, the left-hand side of the above sequent (\ref{eq:trick2}) becomes the following:
$$
\diamondsuit(\bigwedge\limits_{S}
 \a_S\wedge\bigwedge\limits_{S'}\forall p(S'))
$$
and the right-hand side becomes the following:
$$
\diamondsuit(\bigwedge\limits_{S}(\a_S\vee\diamondsuit\bigwedge N(\Box\Sigma,\Sigma;\emptyset))\wedge\bigwedge\limits_{S'}\forall p(S')).
$$

\par\noindent
Observe, that $N(\Box\Sigma|\Sigma;\emptyset)$ is equivalent to $N(\Box\Gamma'|\Gamma';\emptyset)$ by $\Sigma^{\circ}=\Gamma'^{\circ}$. The result now follows by Lemma \ref{lem:k4trick}, putting 
$\be=\bigwedge\limits_{S'}\forall p(S')$ and $\de=\bigwedge N(\Box\Gamma'|\Gamma';\emptyset)$.
\end{rem}

\medskip\par\noindent
\textbf{(i).} The item (i) follows easily by induction on $(\Box\Sigma |\Gamma;\Delta) $ just
because we never add $p$ during the~definition of the~formula
$\forall p(\Box\Sigma|\Gamma;\Delta)$.

\medskip\par\noindent
\textbf{(ii).} we proceed by
induction on the~complexity of $(\Box\Sigma |\Gamma;\Delta)$ given by the~measure function used above to prove termination, and prove that 
$$
\vdash_{G_{Grz}^+}\Box\Sigma|\Gamma,\forall p(\Box\Sigma |\Gamma;\Delta)\Rightarrow\Delta.
$$ 
First let $(\Box\Sigma|\Gamma\Rightarrow\Delta)$ be a~noncritical sequent. Then sequents $(\Box\Sigma_i|\Gamma_{i}\Rightarrow\Delta_{i})\in Cl(\Box\Sigma |\Gamma;\Delta)$ are of lower complexity and by the~induction hypotheses
$$
\vdash_{G_{Grz}^+}\Box\Sigma_i|\Gamma_{i},\forall p(\Box\Sigma_i|\Gamma_{i};\Delta_{i})\Rightarrow\Delta_{i}
$$
for each i. Then by admissibility of weakening and by Lemma \ref{lem:cl+}
$$
\vdash_{G_{Grz}^+}\Sigma|\Gamma,\forall p(\Box\Sigma_1|\Gamma_{1};\Delta_{1}),\ldots,\forall p(\Box\Sigma_k|\Gamma_{k};\Delta_{k})\Rightarrow\Delta,
$$
therefore by a $\wedge$-l inference
$$
\vdash_{G_{Grz}^+}\Box\Sigma|\Gamma,\bigwedge\limits_{(\Box\Sigma_i|\Gamma_{i}\Rightarrow\Delta_{i})\in
Cl(\Box\Sigma|\Gamma;\Delta)}\forall p(\Box\Sigma_i|\Gamma_{i};\Delta_{i})\Rightarrow\Delta,
$$ 
which is by (\ref{eq:defUI2})
$$
\vdash_{G_{Grz}^+}\Gamma,\forall p(\Box\Sigma|\Gamma;\Delta)\Rightarrow\Delta.
$$

\medskip\par\noindent
Let $(\Box\Sigma |\Gamma\Rightarrow\Delta)$ be a~critical sequent
matching the~line 1 of Table \ref{tab:Grz}. Then either (ii) is an~axiom in the
case that $p\in\Pi\cap\Lambda$ or $\bot\in\Pi$, or (ii) is provable
in the case that $\Gamma'\cap\Delta'\neq\emptyset$.

\medskip\par\noindent
Let $(\Box\Sigma |\Gamma\Rightarrow\Delta)$ be a~critical sequent
matching the~line 2 of Table \ref{tab:Grz}. 
We prove 
$$
\vdash_{G_{Grz}^+}\Box\Gamma'|\Pi,\de\Rightarrow\Box\Delta',\Lambda
$$ 
for each disjunct $\de$ used in the line 2 of Table \ref{tab:Grz} to define the interpolant.
\begin{itemize}
\item[-] For each $r\in\Pi$ obviously $\vdash_{G_{Grz}^+}\Box\Gamma'|\Pi,\neg
r,\Box\Gamma'\Rightarrow\Box\Delta',\Lambda$, therefore \\ 
$\vdash_{G_{Grz}^+}\Box\Gamma'|\Pi,\bigvee\limits_{r\in\Pi}\neg
r,\Box\Gamma'\Rightarrow\Box\Delta',\Lambda$.

\item[-] for each $q\in\Lambda$ obviously $\vdash_{G_{Grz}^+}\Box\Gamma'|\Pi
,q,\Box\Gamma'\Rightarrow\Box\Delta',\Lambda$, therefore \\ $\vdash_{G_{Grz}^+}\Box\Gamma'|\Pi,\bigvee\limits_{q\in\Lambda} q,\Box\Gamma'\Rightarrow\Box\Delta',\Lambda$ 

\item[-] For each $\be\in\Delta '$ with $D(\be)\notin\Box\Gamma'$
we have
$$
\vdash_{G_{Grz}^{+}}\Box\Gamma',D(\be)|\Gamma',\forall p(\Box\Gamma',D(\be)|\Gamma';\be)\Rightarrow\be
$$ 
by the~induction hypothesis, which gives
$$
\vdash_{G_{Grz}^{+}}\Box\Gamma ',D(\be),\Box\forall p(\Box\Gamma',D(\be) |\Gamma ';\be)|\Gamma ',\forall p(\Box\Gamma',D(\be) |\Gamma ';\be) \Rightarrow\be
$$ 
by admissible weakening inferences. This yields 
$$
\vdash_{G_{Grz}^{+}}\Box\Gamma',\Box\forall p(\Box\Gamma',D(\be)|\Gamma';\be)|\Pi\Rightarrow \Box\Delta',\Lambda
$$ 
by a~$\Box_{Grz2}^{+}$
inference. Then by weakening and $\Box_{T}^{+}$ inferences
$$
\vdash_{G_{Grz}^{+}}\Box\Gamma '|
\Box\forall p(\Box\Gamma',D(\be)|\Gamma ';\be),\Pi\Rightarrow \Box\Delta',\Lambda.
$$

\item[-] For each $\be\in\Delta '$ with $D(\be)\in\Box\Gamma'$ we have 
$$
\vdash_{G_{Grz}^{+}}\Box\Gamma'|\forall p(\Box\Gamma'|\emptyset;\be)\Rightarrow\beta
$$
by the~induction hypothesis, which gives
$$
\vdash_{G_{Grz}^{+}}\Box\Gamma',\Box\forall p(\Box\Gamma'|\emptyset;\be)|\forall p(\Box\Gamma'|\emptyset;\be),\Gamma'\Rightarrow\beta
$$
by admissible weakening inferences. This yields
$$
\vdash_{G_{Grz}^{+}}\Box\Gamma',\Box\forall p(\Box\Gamma'|\emptyset;\be)|\Pi\Rightarrow \Box\Delta',\Lambda
$$
by a~$\Box_{Grz2}^{+}$ inference and a weakening (notice there is an occurrence of $D(\be)$ missing in $\Gamma'$). Then by weakening and $\Box_{T}^{+}$ inferences
$$
\vdash_{G_{Grz}^{+}}\Box\Gamma '|
\Box\forall p(\Box\Gamma'|\emptyset;\be),\Pi\Rightarrow \Box\Delta',\Lambda.
$$
\end{itemize}
It remains to be proved that 
$$
\vdash_{G_{Grz}^{+}}\Box\Gamma'|\Pi,\diamondsuit \bigwedge N(\Box\Gamma'|\Gamma';\emptyset)\Rightarrow\Box\Delta',\Lambda.
$$
For each $(\Box\Sigma|\Upsilon\Rightarrow\Box\Omega,\Theta)\in
Cl(\Box\Gamma'|\Gamma';\emptyset)$ of the first line of Table \ref{tab:NGrz}, we know by the~induction hypotheses that
$$
\vdash_{G_{Grz}^{+}}\Box\Sigma|\Upsilon,\forall p(\Box\Sigma|\Upsilon;\Box\Omega,\Theta)\Rightarrow\Box\Omega,\Theta.
$$
For each $(\Box\Sigma|\Upsilon\Rightarrow\Box\Omega,\Theta)\in
Cl(\Box\Gamma '|\Gamma';\emptyset)$ of the second line of Table \ref{tab:NGrz}
we have the following:
\begin{itemize}
\item[-]
$
\vdash_{G_{Grz}^{+}}\Box\Sigma|\Upsilon,\bigvee\limits_{q\in\Theta}q\Rightarrow\Box\Omega,\Theta.
$
\item[-] 
$
\vdash_{G_{Grz}^{+}}\Box\Sigma|\Upsilon,\bigvee\limits_{\neg
r\in\Upsilon}\neg r\Rightarrow\Box\Omega,\Theta.
$
\item[-] for each
$\be\in\Omega$ with $D(\be)\notin\Box\Sigma$ by the induction hypotheses
$$
\vdash_{G_{Grz}^{+}}\Box\Sigma,D(\be)|\Sigma,\forall p(\Box\Sigma,D(\be)|\Sigma;\be)\Rightarrow\be
$$ 
and by weakening
and a~$\Box_{Grz2}^+$ inference 
$$
\vdash_{G_{Grz}^{+}}\Box\Sigma,\Box\forall p(\Box\Sigma,D(\be)|\Sigma;\be)|\emptyset\Rightarrow\Box\be.
$$
and by weakening and $\Box_T^+$
$$
\vdash_{G_{Grz}^{+}}\Box\Sigma|\Box\forall p(\Box\Sigma,D(\be)|\Sigma;\be)\Rightarrow\Box\be.
$$

\item[-] for each
$\be\in\Omega$ with $D(\be)\in\Box\Sigma$ by the induction hypotheses
$$
\vdash_{G_{Grz}^{+}}\Box\Sigma|\emptyset,\forall p(\Box\Sigma|\emptyset;\be)\Rightarrow\be
$$ 
and by weakening
and a~$\Box_{Grz2}^+$ inference 
$$
\vdash_{G_{Grz}^{+}}\Box\Sigma,\Box\forall p(\Box\Sigma|\emptyset;\be)|\emptyset\Rightarrow\Box\be.
$$
and by weakening and $\Box_T^+$
$$
\vdash_{G_{Grz}^{+}}\Box\Sigma|\Box\forall p(\Box\Sigma|\emptyset;\be)\Rightarrow\Box\be.
$$
\end{itemize}
Together this yields, using $\vee$-l inferences,
$$
\vdash_{G_{Grz}^{+}}\Box\Sigma|\Upsilon, \bigvee\limits_{q\in\Theta}\vee\bigvee\limits_{\neg r\in\Upsilon}
\vee\bigvee\limits_{D(\be)\notin\Box\Sigma}\Box\forall p(\Box\Sigma,D(\be)|\Sigma;\be)\vee\bigvee\limits_{D(\be)\in\Box\Sigma}\Box\forall p(\Box\Sigma|\emptyset;\be)\Rightarrow\Box\Omega,\Theta.
$$

\medskip\par\noindent
Therefore finally, putting things together for each $(\Box\Sigma|\Upsilon\Rightarrow\Box\Omega,\Theta)\in
Cl(\Box\Gamma '|\Gamma ';\emptyset)$, we obtain, using weakening and $\wedge$-r inferences,
$$
\vdash_{G_{Grz}^{+}}\Box\Sigma,\Upsilon,\bigwedge N(\Box\Gamma'|\Gamma';\emptyset)\Rightarrow\Box\Omega,\Theta.
$$
Now, by closure properties in Lemma \ref{lem:cl}, 
$$
\vdash_{G_{Grz}^{+}}\Box\Gamma '|\Gamma',\bigwedge N(\Box\Gamma'|\Gamma';\emptyset)\Rightarrow\emptyset.
$$ 
By negation and weakening inferences 
$$
\vdash_{G_{Grz}^{+}}\Box\Gamma',D(\neg\bigwedge N(\Box\Gamma',\Gamma';\emptyset))|\Gamma'\Rightarrow\neg\bigwedge N(\Box\Gamma',\Gamma';\emptyset)
$$ 
and by a~$\Box_{Grz2}^+$ inference
$$
 \vdash_{G_{Grz}^{+}}\Box\Gamma'|\Pi\Rightarrow\Box\neg
\bigwedge N(\Box\Gamma',\Gamma';\emptyset),\Box\Delta',\Lambda.
$$
Now, using a~negation inference again, we
obtain 
$$
\vdash_{G_{Grz}^{+}}\Box\Gamma'|\Pi,\diamondsuit\bigwedge N(\Box\Gamma',\Gamma';\emptyset)\Rightarrow\Box\Delta',\Lambda.
$$

\medskip\par\noindent
Putting finally all the above disjuncts together for a critical sequent $(\Box\Gamma'|\Pi;\Box\Delta',\Lambda)$ yields, using $\vee$-l inferences and  the~line 2 of Table \ref{tab:Grz},
$$
\vdash_{G_{Grz}^{+}}\Box\Gamma'|\Pi,\forall p(\Box\Gamma'|\Pi;\Box\Delta',\Lambda)\Rightarrow\Box\Delta',\Lambda.
$$
\textbf{(iii)} We proceed by induction on the~height of the~proof of the 
sequent $(\Box\Theta,\Box\Sigma|\Phi,\Gamma\Rightarrow\Psi,\Delta) $ in
$G_{Grz}^{+}$, and sub-induction on the measure of the sequent $(\Box\Sigma|\Gamma;\Delta)$. We show that 
$$
\vdash_{G_{Grz}^{+}}\Box\Theta|\Phi\Rightarrow\forall p(\Box\Sigma|\Gamma;\Delta),\Psi.
$$
Let us first consider the~last step of the~proof of
$(\Box\Theta,\Box\Sigma|\Phi,\Gamma\Rightarrow\Psi,\Delta)$ is
an~axiom, or, if it is not an axiom, then $(\Box\Sigma|\Gamma;\Delta)$ is a noncritical sequent. In this case we proceed similarly as in Theorem \ref{thm:GL} (iii), the third multiset makes no difference here.

\medskip\par\noindent
Let us then consider that the~last inference of the~proof of
$(\Box\Theta,\Box\Sigma|\Phi,\Gamma\Rightarrow\Psi,\Delta)$ is
a~$\Box_{Grz2}^{+}$ inference. There are two cases to distinguish:
\begin{itemize}
\item[-] Consider first the case when the~principal formula $\Box\a\in\Delta$. Then
the~proof ends with:
\begin{prooftree}
\AxiomC{$\Box\Theta,\Box\Sigma,D(\a)|\Theta,\Sigma\Rightarrow\a$}
\RightLabel{$\Box_{Grz2}^{+}$}
\UnaryInfC{$\Box\Theta,\Box\Sigma|\Gamma,\Phi\Rightarrow\Box\a,\Delta ',\Psi$}
\end{prooftree}
where $\Box\a,\Delta '$ is $\Delta$. Consider $\Box\Sigma\cap\Delta
=\emptyset$ (otherwise $\forall p(\Box\Sigma|\Gamma;\Delta)\equiv\top$ and (iii) holds). Then by the~induction hypotheses
$$
\vdash_{G_{Grz}^{+}}\emptyset|\Box\Theta,\Theta\Rightarrow
\forall p(\Box\Sigma,D(\a)|\Sigma;\a).
$$ 
By invertibility of 
$\Box_{T}^{+}$ inferences, by contraction inferences, and weakening
$$
\vdash_{G_{Grz}^{+}}\Box\Theta,D(\forall p(\Box\Sigma,D(\a)|\Sigma;\a))|\Theta\Rightarrow\forall p(\Box\Sigma,D(\a)|\Sigma;\a),
$$ 
Now, by a~$\Box_{Grz2}^{+}$ inference, we obtain
$$
\vdash_{G_{Grz}^{+}}\Box\Theta|\Phi\Rightarrow\Box\forall p(\Box\Sigma,D(\a)|\Sigma;\a),\Psi
.$$  
By weakening inferences
$$
\vdash_{G_{Grz}^{+}}\Box\Theta|\Theta,\Phi\Rightarrow\Box\forall p(\Box\Sigma,D(\a)|\Sigma;\a),\Psi.
$$
By $\Box_{T}^{+}$ inferences we obtain
$$
\vdash_{G_{Grz}^{+}}\emptyset|\Box\Theta,\Phi\Rightarrow\Box\forall p(\Box\Sigma,D(\a)|\Sigma;\a),\Psi.
$$
By the~line 2 of Table \ref{tab:Grz} and invertibility of the~$\vee$-l rule 
$$
\vdash_{G_{Grz}^{+}}\emptyset|\Box\forall p(\Box\Sigma,D(\a)|\Sigma;\a)\Rightarrow\forall p(\Box\Sigma
|\Gamma;\Box\a,\Delta').
$$ 
The~two sequents above yield (iii) by admissibility of the~cut rule in $G_{Grz}^{+}$.

\medskip
\item[-] Consider next the case when the~principal formula $\Box\a\in\Psi$, i.e.,
$\a$ doesn't contain $p$. Then the~proof ends with:
\begin{prooftree}
\AxiomC{$\Box\Theta,\Box\Sigma,D(\a)|\Theta,\Sigma\Rightarrow\a$}
\RightLabel{$\Box_{Grz2}^{+}$}
\UnaryInfC{$\Box\Theta,\Box\Sigma|\Gamma,\Phi \Rightarrow\Delta,\Box\a,\Psi'$}
\end{prooftree}
where $\Box\a,\Psi'$ is $\Psi$.
Then by the~induction hypotheses
$$
\vdash_{G_{Grz}^{+}}\emptyset|D(\a),\Box\Theta,\Theta\Rightarrow\forall p(\Box\Sigma|\Sigma;\emptyset),\a.
$$
By invertibility of $\Box_{T}^{+}$ inferences and by contraction inferences we obtain
$$
\vdash_{G_{Grz}^{+}}D(\a),\Box\Theta|(\a\rightarrow\Box\a),\Theta\Rightarrow
\forall p(\Box\Sigma|\Sigma;\emptyset),\a.
$$ 
To get rid of
$(\a\rightarrow\Box\a)$, which is $(\neg\a\vee\Box\a)$, we use
invertibility of the $\vee$-l and $\neg$-l rules, and contraction,
to obtain
$$
\vdash_{G_{Grz}^{+}}D(\a),\Box\Theta|\Theta\Rightarrow\forall p(\Box\Sigma|\Sigma;\emptyset),\a.
$$ 
By a~$\neg$-l inference
and weakening
$$
\vdash_{G_{Grz}^{+}}D(\a),\Box\Theta,\Box\neg
\forall p(\Box\Sigma|\Sigma;\emptyset)|\Theta,\neg
\forall p(\Box\Sigma|\Sigma;\emptyset)\Rightarrow\a.
$$ 
By
a~$\Box_{Grz2}^{+}$ inference
$$
\vdash_{G_{Grz}^{+}}\Box\Theta,\Box\neg
\forall p(\Box\Sigma|\Sigma;\emptyset)|\Phi\Rightarrow\Box\a,\Psi'.
$$ 
Since weakening is admissible in $Gm_{Grz}^{+}$, we obtain
$$
\vdash_{G_{Grz}^{+}}\Box\Theta,\Box\neg
\forall p(\Box\Sigma|\Sigma;\emptyset)|\Theta,\neg
\forall p(\Box\Sigma|\Sigma;\emptyset),\Phi \Rightarrow\Box \a,\Psi'
$$ 
and now $\Box_{T}^{+}$ inferences and a~$\neg$-l inference yield
$$
\vdash_{G_{Grz}^{+}}\emptyset|\Box\Theta,\Phi
\Rightarrow\diamondsuit\forall p(\Box\Sigma|\Sigma;\emptyset),\Box
\a,\Psi'.
$$ 
By weakening inferences
$$
\vdash_{G_{Grz}^{+}}\Box\Theta|\Theta,\Phi\Rightarrow\diamondsuit
\forall p(\Box\Sigma|\Sigma;\emptyset),\Box\a,\Psi'.
$$ By $\Box_{T}^{+}$ inferences
$$
\vdash_{G_{Grz}^{+}}\emptyset|\Box\Theta,\Phi \Rightarrow\diamondsuit
\forall p(\Box\Sigma|\Sigma;\emptyset),\Box\a,\Psi'.
$$ 
By
the~line 2 of Table \ref{tab:Grz}, invertibility of the~$\vee$-l rule, and by (\ref{eq:trick2}) we have
$$
\vdash_{G_{Grz}^{+}}\emptyset|\diamondsuit\forall p(\Box\Sigma|\Sigma;\emptyset)\Rightarrow\forall p(\Box\Sigma|\Sigma;\emptyset).
$$ 
The~two sequents above yield (iii) by
admissibility of the~cut rule in $G_{Grz}^{+}$.
\end{itemize}
Let us consider that the~last inference of the~proof of
$(\Box\Theta,\Box\Sigma|\Phi,\Gamma\Rightarrow\Psi,\Delta) $ is
a~$\Box_{Grz1}^{+}$ inference. Again, we distinguish two cases:
\begin{itemize}
\item[-] Consider first the case when the~principal formula $\Box\a\in\Delta$. Then
the~proof ends with:
\begin{prooftree}
\AxiomC{$\Box\Theta,\Box\Sigma|\emptyset\Rightarrow\a$}
\RightLabel{$\Box_{Grz1}^{+}$}
\UnaryInfC{$\Box\Theta,\Box\Sigma|\Gamma,\Phi\Rightarrow\Box\a,\Delta',\Psi$}
\end{prooftree}
where $\Box\a,\Delta '$ is $\Delta$. Consider $\Box\Sigma\cap\Delta
=\emptyset$ (otherwise $\forall p(\Box\Sigma|\Gamma;\Delta)\equiv\top$ and (iii)
holds).
Then by the~induction hypotheses
$$
\vdash_{G_{Grz}^{+}}\emptyset|\Box\Theta\Rightarrow
\forall p(\Box\Sigma|\emptyset;\a).
$$ 
By invertibility of $\Box_{T}^{+}$ inferences,
by contraction inferences, and weakening
$$
\vdash_{G_{Grz}^{+}}\Box\Theta,D(\forall p(\Box\Sigma|\emptyset;\a))|\Theta\Rightarrow
\forall p(\Box\Sigma|\emptyset;\a),
$$ 
Now, by a~$\Box_{Grz2}^{+}$ inference,
we obtain
$$
\vdash_{G_{Grz}^{+}}\Box\Theta|\Phi\Rightarrow\Box\forall p(\Box\Sigma|\emptyset;\a),\Psi.
$$  
By weakening inferences
and $\Box_{T}^{+}$ inferences
we obtain
$$
\vdash_{G_{Grz}^{+}}\emptyset|\Box\Theta,\Phi \Rightarrow\Box
\forall p(\Box\Sigma|\emptyset;\a)\Psi.
$$ 
By the~line 2 of Table \ref{tab:Grz} and
invertibility of the~$\vee$-l rule we have
$$
\vdash_{G_{Grz}^{+}}\emptyset|\Box\forall p(\Box\Sigma|\emptyset;\a)\Rightarrow \forall p(\Box\Sigma|\Gamma;\Box\a,\Delta') .
$$ 
The~two sequents above yield (iii) by
admissibility of the~cut rule in $G_{Grz}^{+}$.

\medskip
\item[-] Consider next the case that the~principal formula $\Box\a\in\Psi$, i.e.,
$\a$ doesn't contain $p$. Then the~proof ends with:
\begin{prooftree}
\AxiomC{$\Box\Theta,\Box\Sigma|\emptyset\Rightarrow\a$}
\RightLabel{$\Box_{Grz1}^{+}$}
\UnaryInfC{$\Box\Theta,\Box\Sigma|\Gamma,\Phi\Rightarrow\Delta,\Box\a,\Psi'$}
\end{prooftree}
where $\Box\a,\Psi'$ is $\Psi$. Then by the~induction hypotheses
$$
\vdash_{G_{Grz}^{+}}\emptyset|\Box\Theta
\Rightarrow \forall p(\Box\Sigma|\emptyset;\emptyset),\a.
$$ 
Notice that
$(\Box\Sigma|\emptyset;\emptyset)$ is a~critical sequent with all but
one multisets empty, and by the~table $\forall p(\Box\Sigma|\emptyset;\emptyset)\equiv\bot$. Thus we have in fact
$$
\vdash_{G_{Grz}^{+}}\emptyset|\Box\Theta\Rightarrow\a.
$$ 
By invertibility of $\Box_{T}^{+}$ inferences we obtain
$$
\vdash_{G_{Grz}^{+}}\Box\Theta|\Theta\Rightarrow\a
$$ 
and by weakening
$$
\vdash_{G_{Grz}^{+}}\Box\Theta,D(\a)| \Theta\Rightarrow\a.
$$ 
By
a~$\Box_{Grz2}^{+}$ inference
$$
\vdash_{G_{Grz}^{+}}\Box\Theta|\Phi\Rightarrow\Box\a,\Psi'.
$$ 
By admissibility of weakening we
obtain 
$$
\vdash_{G_{Grz}^{+}}\Box\Theta|\Phi\Rightarrow\forall p(\Box\Sigma|\Gamma;\Delta),\Box\a,\Psi'.
$$
\end{itemize}
\end{proofof}

\subsection{Concluding remarks} We have provided an effective construction of uniform interpolants in provability logics. We would like to point out, that even if the proofs as presented are not fully constructive, the only part that is not constructive is the completeness of the two calculi without the cut rule. This can be, in both cases, repaired by completing the proof search argument and make it into a decision procedure. 

What we also left open in this paper is to investigate which distribution laws the quantifiers satisfy. For example, it is the case that, in the basic modal logic \textbf{K}, the universal bisimulation quantifier commutes with the diamond modality \cite{Bil06},\cite{Fre06}. In fact, it commutes with (the dual of) the cover modality, which is a principle that, besides the usual axioms and rules for quantification, axiomatizes bisimulation quantifiers over \textbf{K}. Whether a similar  insight can be obtained for \textbf{GL} is not clear at the moment.

\bibliographystyle{amsplain}
\bibliography{UIbib}
\end{document}